\documentclass[12pt]{iopart}
\usepackage{iopams,setstack}
\usepackage{stmaryrd}
\usepackage{amssymb,amsthm}
\usepackage{graphicx,float}
\usepackage[colorlinks=true,linkcolor=blue,citecolor=blue,urlcolor=blue]{hyperref}
\usepackage{color} 
\usepackage{mathrsfs}
\usepackage{subfigure}
\newtheorem{theorem}{Theorem}[section]

\makeatletter
\newcommand{\bigintss}{\@ifnextchar_\@bigintsssub\@bigintssnosub}
\def\@bigintsssub_#1{\def\@int@subscript{#1}\@ifnextchar^\@bigintsssubsup\@bigintsssubnosup}
\def\@bigintsssubsup^#1{\mathop{\text{\LARGE$\int_{\text{\normalsize$\scriptstyle\kern-0.25em\@int@subscript$}}^{\text{\normalsize$\scriptstyle#1$}}$}}\nolimits}
\def\@bigintsssubnosup{\mathop{\text{\LARGE$\int_{\text{\normalsize$\scriptstyle\@int@subscript$}}$}}\nolimits}
\def\@bigintssnosub{\@ifnextchar^\@bigintssnosubsup\@bigintssnosubnosup}
\def\@bigintssnosubsup^#1{\mathop{\text{\LARGE$\int^{\text{\normalsize$\scriptstyle#1$}}$}}\nolimits}
\def\@bigintssnosubnosup{\mathop{\text{\LARGE$\int$}}\nolimits}
\makeatother

\begin{document}
\title{A physiology--based parametric imaging method for FDG--PET data.}
\author{Mara Scussolini$^1$, Sara Garbarino$^2$, Gianmario Sambuceti$^3$, Giacomo Caviglia$^1$ and Michele Piana$^4$}
\address{$^1$ Dipartimento di Matematica, Universit\`a di Genova, Genova, Italy}
\address{$^2$ Centre for Medical Image Computing, Department of Computer Science, University College London, UK}
\address{$^3$ Dipartimento di Medicina Nucleare, IRCCS-IST San Martino, and Dipartimento di Scienze della Salute, Universit\`a di Genova, Genova, Italy.}
\address{$^4$ Dipartimento di Matematica, Universit\`a di Genova, and CNR-SPIN, Genova, Italy.}

\ead{\mailto{scussolini@dima.unige.it}}

\noindent{\it Keywords\/}: numerical inverse problems, parametric imaging, compartmental analysis, nuclear medicine data 

\begin{abstract}
Parametric imaging is a compartmental approach that processes nuclear imaging data to estimate the spatial distribution of the kinetic parameters governing tracer flow. The present paper proposes a novel and efficient computational method for parametric imaging which is potentially applicable to several compartmental models of diverse complexity and which is effective in the determination of the parametric maps of all kinetic coefficients. We consider applications to [$^{18}$F]-fluorodeoxyglucose Positron Emission Tomography (FDG--PET) data and analyze the two--compartment catenary model describing the standard FDG metabolization by an homogeneous tissue and the three--compartment non--catenary model representing the renal physiology. We show uniqueness theorems for both models. The proposed imaging method starts from the reconstructed FDG--PET images of tracer concentration and preliminarily applies image processing algorithms for noise reduction and image segmentation. The optimization procedure solves pixel--wise the non--linear inverse problem of determining the kinetic parameters from dynamic concentration data through a regularized Gauss--Newton iterative algorithm. The reliability of the method is validated against synthetic data, for the two--compartment system, and experimental real data of murine models, for the renal three--compartment system. 
\end{abstract}

\section{Introduction}
Nuclear Medicine analyzes dynamic data of functional processes related to a specific metabolic activity. Nuclear medicine imaging data are acquired by means of devices that detect the product of the decay of radioisotopes in a radioactive tracer, which is bound to molecules with known biological properties and diffused in the living organism. Positron Emission Tomography (PET) \cite{Bailey,Ollinger} is the most modern nuclear medicine modality and FDG--PET is the PET modality in which the radiopharmaceutical [$^{18}$F]-fluorodeoxyglucose (FDG) is used as a tracer to evaluate glucose metabolism \cite{Kelloff}. FDG--PET dynamic images of tracer distribution, obtained from the measured radioactivity by applying an appropriate reconstruction algorithm, are a reliable estimate of the functional behavior of glucose into tissues. Therefore, FDG--PET experiments allow clinicians to detect and stage diseases related to the pathological glucose consumption, such as cancer \cite{Annibaldi,Cairns,Warburg} or diabetes \cite{Basu,Iozzo}. 

Compartmental analysis \cite{Gunn,Schmidt} is a powerful tool for processing dynamic PET data and estimating physiological kinetic parameters explaining the tracer metabolism. This kind of analysis relies on (I) the construction of a forward model (typically in the form of a Cauchy problem) for tracer concentration, in which the differential equations' coefficients are the kinetic parameters, and on (II) the application of an inversion technique to retrieve such kinetic parameters from measurements of the dynamic tracer concentration. 

Compartmental analysis can be mainly subdivided in two classes: Region Of Interest (ROI) kinetic modeling and parametric imaging. ROI--based methods \cite{Carson,Garbarino_kidney,Garbarino_liver,Vanzi} return a single set of kinetic parameters for a homogeneous region of tissue, whose Time Activity Curve (TAC) is obtained averaging the PET activity over the region at each time frame. On the other hand, parametric imaging \cite{Reader} aims at evaluating the set of kinetic parameters for every pixel of the PET images, thus providing the spatial distribution of each model parameter. This approach is particularly useful when the tissue under examination cannot be effectively segmented into homogeneous regions that could be modeled with a single kinetic parameter set. There exist \emph{indirect} and \emph{direct} approaches to parametric imaging. Indirect methods work by first reconstructing the dynamic PET images and then estimating the kinetic parameters at each pixel \cite{Zhou1,Zhou2}. Direct methods estimate directly the space-varying kinetic parameters from the measured PET sinogram \cite{WangQi}. The direct approach has been proved to reduce the Signal to Noise Ratio (SNR) with respect to indirect techniques \cite{OSullivan}, although strongly relies on the implementation of an efficient inversion algorithm capable of reconstructing the parameters on a dense set of pixels \cite{Kamasak}. However, most parametric imaging methods (both direct and indirect) rely on linearized compartmental models and/or provide parametric images of algebraic combinations of the kinetic coefficients \cite{Logan1,Logan2,Patlak,Thie,Tsoumpas}. Rather few methods are able to reconstruct maps of each single parameter, and most of them consider simple one-- and two--compartment models \cite{Carson_parametric,Huesman,Kamasak,Limber}. 

 In this paper, we want to exploit recent advances in ROI--based analysis \cite{Garbarino_kidney,Garbarino_liver,Vanzi,Watabe} in order to realize a novel, computationally efficient imaging procedure that can be used for parametric imaging in the case of complex physiological systems. The novelties of our approach can be summarized as follows.
\begin{itemize}
\item The method can address non--standard multi--compartment physiologies on a large set of pixels, can do this by exploiting complex physiological information and can realize the analysis in a parametric imaging framework. At this stage, the overall approach requires the use of reconstructed dynamic images but, differently than most (if not all) direct methods, which consider just one ore two catenary compartments, here more complex physiological systems can be considered.
\item From a mathematical perspective, specific physiological conditions assure the identifiability of the compartmental inverse problem. To the best of our knowledge, the uniqueness result given in this paper for a non--catenary three--compartment model is the most general one currently available in compartmental analysis.
\item From a technological viewpoint, the paper introduces an overall automatic pipeline that comprises both pre-processing and inversion, whereby automation has been realized also thanks to an implementation of the Gauss--Newton algorithm that works pixel--wise, with an optimized choice of the regularization parameter based on Generalized Cross Validation. Such pipeline takes as input the dynamic PET images, gives the user the possibility to select the organ and the related physiological compartmental model, realizes a standard but reliable image pre-processing for noise-reduction and image segmentation, and finally computes  the values of a notable number of kinetic parameters, pixel--wise.
\end{itemize}
What we obtain is a method general enough to work for both two- and three--compartment models, effective enough to provide maps of all the kinetic coefficients involved, and that therefore can be in principle extended to envisage more than one model for a single PET image and physiologies described by more than three compartments. The proposed pipeline is validated on synthetic data mimicking a standard two--compartment system for a generic homogeneous tissue and is applied against experimental measurements concerning the renal system of mouse models. For this latest case, we also provide, for the first time, proof of identifiability valid for non--catenary models, as the one describing the renal system.

The scheme of the paper is as follows. In Section 2 the two--compartment and three--compartment models are presented and the related forward problems are described. Section 3 deals with the formalization of the respective inverse problems and with the discussion of identifiability of the models. Section 4 describes the computational parametric imaging method. Section 5 provides the numerical validation of the computational method in the case of the two--compartment catenary model and then applies the method against experimental murine data for the analysis of the three--compartment non--catenary model of the renal system. Our conclusions are offered in Section 6.

\section{Mathematical models} \label{sec:models}
Compartmental analysis \cite{Gunn,Schmidt} identifies different functional compartments in the physiological system of interest, each one associated with a specific metabolic state of the tracer.
The tracer typically is injected into the blood and the tracer concentration in the blood is mathematically modeled by the Input Function (IF) of the compartmental system.  In this work, the IF is assumed to be known as it can be obtained by drawing ROIs on
reconstructed PET images at different time points, in correspondence with the left ventricle. When the IF is not given, suitable reference tissue models have to be taken into account \cite{Lammertsma,Schmidt,Zanotti,Zhou3}. The time dependent concentrations of tracer in each compartment constitute the state variables of the model and can be determined from PET data. The time evolution of the state variables, i.e. the kinetics of the system, is here modeled by a linear system of Ordinary Differential Equations (ODEs) with constant coefficients, expressing the conservation of tracer during the flow between compartments. The coefficients define the input/output rate of tracer for each compartment and represent the physiological parameters describing the metabolism of the system. The forward problem is constructed by assuming that the tracer coefficients are known and by solving the system of ODEs for the unknown concentrations. We are aware that recent works \cite{PDE} take into account macroscopic flow conditions (especially to model cardiac perfusion), introducing a PDE--based framework. However, the simplifying assumptions of time independence of parameters and of no spatial exchange between compartments are well established for FDG--PET analysis \cite{Wernick}, and are the conditions assumed to hold throughout this paper.

This section is devoted to the description of the standard two--compartment catenary model and a three--compartment non--catenary model developed for the renal physiology and to the discussion of the related forward problems. In the following analysis, we denote with capital $C$ (kBq ml$^{-1}$) the concentrations and with the related suffixes the corresponding compartment; we use the notation $k_{yx}$ (min$^{-1}$) for the kinetic parameter describing the tracer exchange to the target compartment $y$ from the source compartment $x$. The kinetic parameters (also known as rate constants or exchange coefficients) are real positive numbers and the plus or minus signs against them characterize incoming and outgoing fluxes, respectively. Notice that the spatial dependence on the pixel index $(i,j)$ in the compartment concentrations and in the kinetic parameters is omitted but implied. We indicate a generic ROI compartment concentration with the apex ROI.

\subsection{Two--compartment catenary system} \label{subsec:2C}

\begin{figure}[htb]  
\centering
\includegraphics[width=10cm]{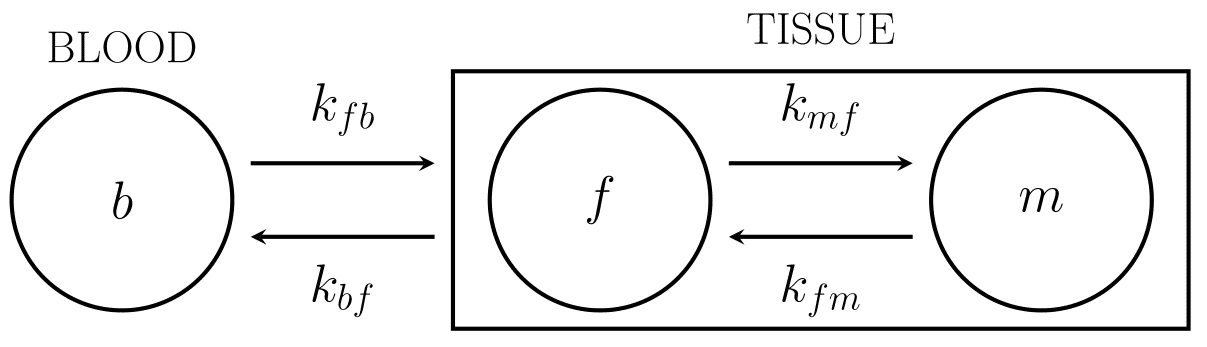}
\caption{The compartmental model for the two--compartment catenary system describing the FDG metabolization in a generic tissue.}
\label{fig:model_2C}
\end{figure}

The compartmental model describing the FDG metabolism of phosphorylation--de-phosphorylation is the two--compartment catenary model shown in \figurename~\ref{fig:model_2C} \cite{Delbary_id,Sokoloff}.

The two--compartment catenary model consists of: 
\begin{itemize} 
\item the blood compartment $b$;
\item two functional compartments: compartment $f$, accounting for free FDG, and compartment $m$, accounting for metabolized FDG;
\item four exchange coefficients: $k_{fb}$ and $k_{bf}$ between the blood and the free compartment, $k_{mf}$ and $k_{fm}$ between the free and the metabolized ones.
\end{itemize}

Balance of tracer concentrations leads to the following system of ODEs:
\begin{equation} \label{eqn:2C_ode} 
\frac{d\boldsymbol{C}}{dt}(t) = \dot{\boldsymbol{C}}(t) = \boldsymbol{M} \boldsymbol{C}(t) + \boldsymbol{w}(t)
\end{equation}
where
\begin{eqnarray*}
\boldsymbol{C} &= \pmatrix{ C_f \cr C_m} , \ \boldsymbol{M} = \pmatrix{-(k_{bf} + k_{mf}) & k_{fm} \cr k_{mf} & -k_{fm}} , \\  \boldsymbol{w} &= k_{fb} C_b^{\text{ROI}} \boldsymbol{e_1} = \pmatrix{ k_{fb} C_b^{\text{ROI}} \cr 0 } , \mbox{ and } \boldsymbol{e_1} = \pmatrix{ 1 \cr 0 }.
\end{eqnarray*}
The initial condition is $C_f(0)=C_m(0)=0$, i.e. $\boldsymbol{C}(0) = 0$. The blood ROI compartment concentration $C_b^{\text{ROI}}$ plays the role of the given IF of the two--compartment model.

The analytical solution $\boldsymbol{C}(t)$ of (\ref{eqn:2C_ode}), formally expressing the forward problem of evaluating the concentrations from the kinetic parameters $\boldsymbol{k} = (k_{fb}, k_{bf}, k_{mf}, k_{fm})^T$, is given by
\begin{equation} \label{eqn:2C_sol}
\boldsymbol{C}(t) = \int_0^t e^{(t-\tau)\boldsymbol{M}} \boldsymbol{w}(\tau) \, d\tau = k_{fb} \int_0^t e^{(t-\tau)\boldsymbol{M}} C_b^{\text{ROI}}(\tau) \boldsymbol{e_1} \, d\tau \ , 
\end{equation}
with the time variable $t \in \mathbb{R}_+$.

\subsection{Three--compartment non--catenary system} \label{subsec:3C}

\begin{figure}[htb]  
\centering
\includegraphics[width=10cm]{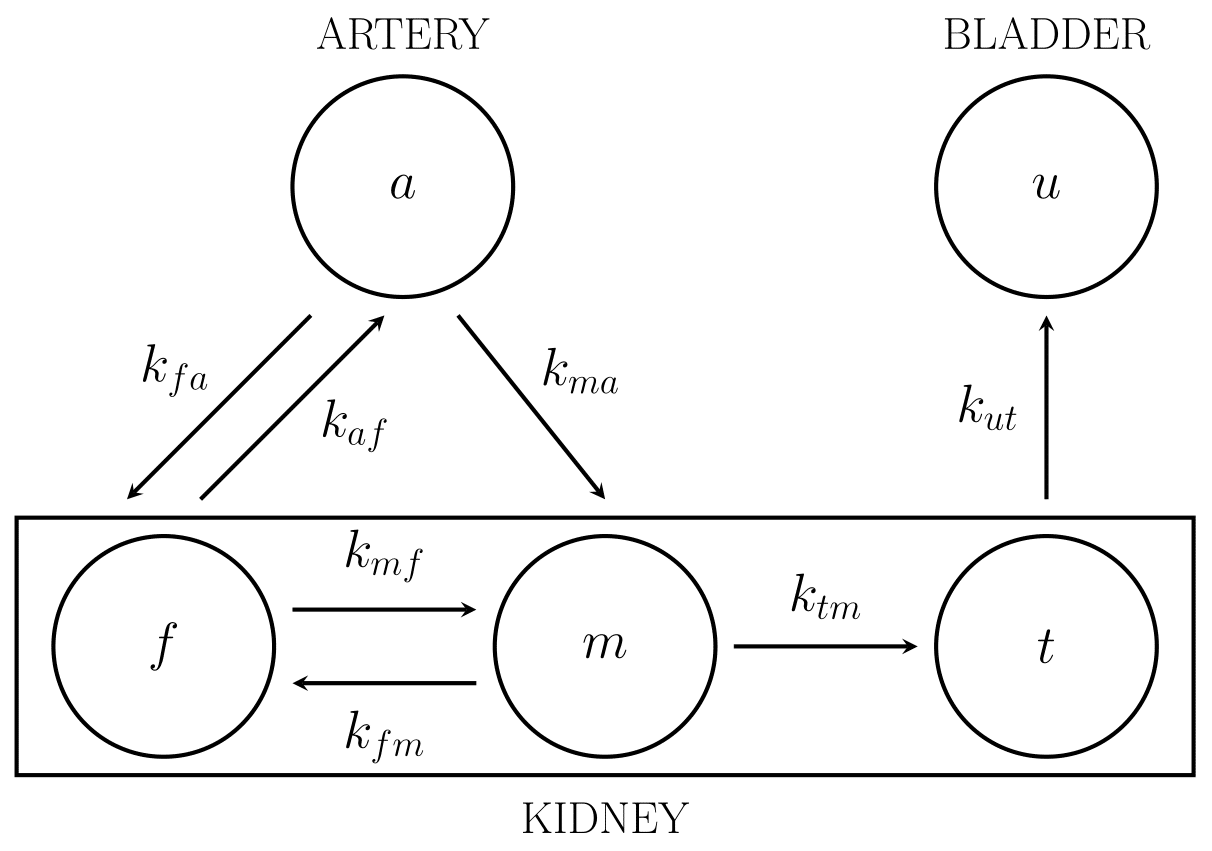}
\caption{The compartmental model for the three--compartment non--catenary system describing the FDG kinetics inside the kidney.}
\label{fig:model_kidney}
\end{figure}

Once injected into the system, the tracer reaches the kidneys and infuses the organs via the blood stream through the renal artery ($a$), in \figurename~\ref{fig:model_kidney}. Here, we consider the usual two--compartment model describing the FDG phosphorylation--de-phosphorylation processes, obtaining the free tracer ($f$) and the metabolized tracer ($m$), both located in the extravascular kidney tissue. However, in order to study the role of the mechanisms transporting glucose back to the metabolism, we need to include the reabsorption compartment ($t$), anatomically identified with the renal tubule. 
Moreover, we add the external urine compartment $u$, anatomically identified with the bladder, accounting for the tracer there accumulated, thanks to the excretion mechanism (differently from glucose, FDG is poorly absorbed in the renal tubule and is largely excreted in the urine, with accumulation in the bladder \cite{Shreve}). 

The resulting three--compartment non--catenary model represented in \figurename~\ref{fig:model_kidney} \cite{Garbarino_kidney,Qiao} has the following kinetic parameters:
\begin{itemize}
\item $k_{fa}$ and $k_{af}$ between the arterial IF and the free FDG compartment;
\item $k_{ma}$ from the arterial compartment to the metabolized FDG one (filtration process);
\item $k_{fm}$ and $k_{mf}$ between the free FDG and the metabolized FDG compartments (phosphorylation and de-phosphorylation processes);
\item $k_{tm}$ from the metabolized FDG compartment to the tubule (reabsorption process);
\item $k_{ut}$ from the tubule compartment to the bladder pool (excretion process).
\end{itemize}

The resulting system of ODEs is
\begin{equation} \label{eqn:kidney_ode} 
\cases{\dot{\boldsymbol{C}}(t) = \boldsymbol{M} \boldsymbol{C}(t) + \boldsymbol{w}(t) \\
\boldsymbol{C}(0) = 0} \ ,
\end{equation}
where
\begin{eqnarray*}
\boldsymbol{C} &= \pmatrix{ C_f \cr C_m \cr C_t } , \ \boldsymbol{M} = \pmatrix{ -(k_{af} + k_{mf}) & k_{fm} & 0 \cr k_{mf} & -(k_{fm} + k_{tm}) & 0 \cr 0 & k_{tm} & -k_{ut} } , \\ 
\boldsymbol{w} &= k_{fa} C_a^{\text{ROI}} \boldsymbol{e_1} + k_{ma} C_a^{\text{ROI}} \boldsymbol{e_2} = \pmatrix{ k_{fa} C_a^{\text{ROI}} \cr k_{ma} C_a^{\text{ROI}} \cr 0 } , 
\end{eqnarray*}
\begin{eqnarray*}
\mbox{ and } \boldsymbol{e_1} = \pmatrix{1 \cr 0 \cr 0 } , \mbox{ and } \boldsymbol{e_2} = \pmatrix{0 \cr 1 \cr 0 }.
\end{eqnarray*}

The analytical solution for $\boldsymbol{C}$, in terms of the kinetic coefficients $\boldsymbol{k} = (k_{fa}, k_{ma}, k_{af}, k_{mf}, k_{fm}, k_{tm}, k_{ut})^T$, is given by
\begin{eqnarray} \label{eqn:kidney_sol1} 
& \boldsymbol{C}(t) = \int_0^t e^{(t-\tau)\boldsymbol{M}} \boldsymbol{w}(\tau) \, d\tau = \\
\label{eqn:kidney_sol2} &= k_{fa} \int_0^t e^{(t-\tau)\boldsymbol{M}} C_a^{\text{ROI}}(\tau) \boldsymbol{e_1} \, d\tau + k_{ma} \int_0^t e^{(t-\tau)\boldsymbol{M}} C_a^{\text{ROI}}(\tau) \boldsymbol{e_2} \, d\tau \ , 
\end{eqnarray}
 with the time variable $t \in \mathbb{R}_+$ and the arterial IF $C_a^{\text{ROI}}$ a known function.
 
\subsubsection{Conditions from physiology.} \label{subsubsec:ktm_kut} 
We expect that in general the three--compartment system described in \ref{subsec:3C} is not identifiable in that it is a non--catenary model involving seven unknown rate constants. Identifiability is achieved by reducing the number of unknowns through the introduction of constraints coming from renal physiology. It is shown in the course of the discussion that the constraints are used effectively in the discussion of uniqueness but they cannot be applied in the procedure for the solution of the inverse problem. 
 
To introduce the restrictions, we first consider tracer kinetics within the ROI--based framework, by looking at the overall bulk flow of tracer between the ROI compartments tubule and bladder. The pertinent properties are then extended to the pixel framework, in analogy with the extension of ROI based features of other compartments, e.g. the dynamical system of tracer exchange between model compartments.  
Specifically, we have to reconsider the role of the bladder, which accumulates urine and tracer in the course of the experiment. The bladder is connected to the kidneys but is not part of the kidneys, and is related to a strictly global behavior of the renal system. The influence of bladder in tracer kinetics has to be inserted into the description of the characteristic features of the pixel constituents of the renal compartment system, although pixels are not physically connected to the bladder. The following points describe a few aspects of bladder involvement in tracer kinetics that are relevant to parametric modeling and are applied in the subsequent analysis of identifiability. We assume that the volume and the activity of bladder are known, although we do not need the corresponding explicit values.
\begin{enumerate}
\item The bladder compartment of tracer density $C_u^{\text{ROI}}$ is the only compartment whose volume significally changes (specifically increases) in time. The condition of tracer balance for bladder is expressed by equation
\begin{equation} \label{eqn:bladder_ode}
\frac{d}{dt} (V_u C_u^{\text{ROI}}) = F^{\text{ROI}}_{ut} C_t^{\text{ROI}} \quad \mbox{with} \quad C_u^{\text{ROI}}(0) = 0 \ ,
\end{equation} 
where $V_u$ (ml) indicates the bladder volume, and $F^{\text{ROI}}_{ut}$ (ml min$^{-1}$) is the bulk flow entering the bladder from the tubule region. In words, equation (\ref{eqn:bladder_ode}) states that the (positive) time rate of the tracer content of bladder equals the amount of tracer carried inward by the flux of urine. 
\item  We assume that the accumulation rate of urine is constant, consistently with the  assumption of stationarity and the condition of resting state of the subject during PET acquisition. Therefore, the growth of the bladder volume $V_u$ during time is linear and hence the flux rate of urine into bladder satisfies \cite{Garbarino_kidney}
\begin{equation} \label{eq:F^R_ut}
F^{\text{ROI}}_{ut} = \frac{V_u(t_f) - V_u(\bar{t})}{t_f - \bar{t}} \ , 
\end{equation} 
where $t_f$ is the end time point and $\bar{t}$ is a generic time instant. In principle, an estimate of $F^{\text{ROI}}_{ut}$ is obtained from experimental values of $V_u$ at different time points. As a consequence we find that
\begin{equation} \label{eqn:bladder_ode_C_t}
C_t^{\text{ROI}}= \frac 1 {F^{\text{ROI}}_{ut}} \frac{d}{dt} (V_u C_u^{\text{ROI}}) \ ,
\end{equation}
showing that $C_t^{\text{ROI}}$ can be estimated from measurements performed on bladder. 
\item Consider the renal tubule ROI compartment. We regard as a result from physiology the existence of a direct proportionality between the bulk ingoing flow $F^{\text{ROI}}_{tm}$ (ml min$^{-1}$), entering  the renal tubule system from the metabolized compartment, and the bulk outgoing flow $F^{\text{ROI}}_{ut}$, directed towards bladder. Accordingly, we can write
\begin{equation} \label{eqn:fluxes}
F^{\text{ROI}}_{tm} = \gamma F^{\text{ROI}}_{ut} \ ,
\end{equation} 
where $\gamma$ is a constant factor establishing the order of magnitude of the bulk quantity passing from the tubule to the bladder. For example, for the mouse model the value of $\gamma$ is typically equal to $10^2$ \cite{kid_phys}. Substitution of (\ref{eq:F^R_ut}) into (\ref{eqn:fluxes}) provides $F^{\text{ROI}}_{tm}$ in terms of data.
\item The balance of tracer flow inside the overall tubule system is given by \cite{Garbarino_kidney}
\begin{equation} \label{eqn:tubule_ode}
\frac{d}{dt}(V_K V_t C_t^{\text{ROI}}) = F^{\text{ROI}}_{tm}C_m^{\text{ROI}} - F^{\text{ROI}}_{ut} C_t^{\text{ROI}} \ ,
\end{equation} 
where $V_K$ (ml) is the renal volume, $V_t$ the fraction of the tubule volume, and thus $V_K V_t$ (ml) is the total volume of the tubule. Equation (\ref{eqn:tubule_ode}) may be written in the equivalent form
\begin{equation} \label{eqn:tubule_ode_bis}
\dot{C}_t^{\text{ROI}} = k^{\text{ROI}}_{tm}C_m^{\text{ROI}} - k^{\text{ROI}}_{ut} C_t^{\text{ROI}} \ ,
\end{equation} 
where the ROI parameters $k^{\text{ROI}}_{tm}$ and $k^{\text{ROI}}_{ut}$ are defined by
\begin{equation} \label{eqn:k_tubule}
k^{\text{ROI}}_{tm} = \frac{F^{\text{ROI}}_{tm}}{V_K V_t} ,\qquad k^{\text{ROI}}_{ut} = \frac{F^{\text{ROI}}_{ut}}{V_K V_t} ,   
\end{equation} 
with
\begin{equation} \label{eqn:k_tubule_bis}
k^{\text{ROI}}_{tm}=\gamma k^{\text{ROI}}_{ut} \ .
\end{equation} 
The last equality follows from equation (\ref{eqn:fluxes}). Also $k^{\text{ROI}}_{tm}$ and $k^{\text{ROI}}_{ut}$ are determined by data through $F^{\text{ROI}}_{tm}$, $F^{\text{ROI}}_{ut}$, $V_K$, and $V_t$. 
\end{enumerate}

We now come to the parametric formulation. As already observed, we reproduce at the pixel level a few relevant features of the ROI approach; in particular, each pixel is regarded as capable interchanging carrier fluid and tracer with bladder. In line with equation  (\ref{eqn:bladder_ode_C_t}), we assume that, at each pixel, 
\begin{equation} \label{eqn:bladder_ode_C_t_bis}
C_t= \frac 1 {F_{ut}} \frac{d}{dt} (V_u C_u^{\text{ROI}}) \ ,
\end{equation}
where $F_{ut}$ describes the contribution to bladder volume change per unit time arising from the flux of fluid through the single pixel. The coefficient $F_{ut}$ is dependent on the position within the renal tissue and is considered as fixed, in analogy with (\ref{eqn:bladder_ode_C_t}).
The correspondent of equation (\ref{eqn:tubule_ode_bis}) is postulated as  
\begin{equation} \label{eqn:kidney_ode_3}
\dot{C}_t = k_{tm} C_m - k_{ut} C_t \ ,
\end{equation} 
accounting for the dynamic of the pixel tubule compartment. As expected, equation (\ref{eqn:kidney_ode_3}) coincides with the third differential equation of the system of ODEs (\ref{eqn:kidney_ode}). 
The rate coefficient $k_{tm}$ and $k_{ut}$ are position dependent, as in (\ref{eqn:kidney_ode}), and are regarded as fixed.

To summarize, for the pixel dependent tubule concentration and rate coefficients, we have obtained from physiology the conditions that $C_t(t)$, $k_{tm}$, $k_{ut}$ are fixed, and that $k_{tm} = \gamma k_{ut}$. 
Such conditions are applied in the proof of uniqueness. It is important to remark that, although $C_t$, $k_{tm}$, and $k_{ut}$ are fixed by physiology, the corresponding time dependence and values are not known explicitly. This means that they are still to be determined through optimization.

\section{Inverse problem: data and solution} \label{sec:inverse}

The forward model equation, as the analytical solution to the system of ODEs, describes the time behavior of the tracer concentration given the TAC for tracer concentration in blood and the exchange coefficients. Given such equation, compartmental analysis requires the determination of the parameters by utilizing measurements of the tracer concentrations provided by PET imaging and applying an optimization scheme for the solution of the inverse problem. The reconstructed PET images of tracer concentration supply information on the IF and provide an estimate of a weighted sum of the concentrations of the different compartments considered, at each time point of the acquisition. 
In this work, we make use of the Gauss--Newton method supplied with a regularizing term \cite{GN, Delbary_meth, Nocedal, Vogel}, to estimate the exchange parameters. It has been shown in \cite{Delbary_meth} that Gauss--Newton regularization in the compartmental framework provides reconstructions of the kinetic parameters that are more stable with respect to Levenberg--Marquardt method (see Tables 1-3 in that paper). Further, the matrix differentiation step required at some stage of the optimization analysis is in this method performed analytically, thus avoiding time consuming numerical differentiation, and tackling one standard drawback of Newton methods. Also, by searching for zeroes of non--linear functionals, Newton methods do not need to \emph{a priori} select a topology in the data space, as in the case of least--squares approaches. On the other hand, this class of methods, as compared to standard optimization-through-minimization techniques, may lack in convergence if the starting point is taken too far away from the ground truth. In the cases we consider, thanks to the identifiability results, we can have good expectations of robustness and accuracy even if the initial guess is not carefully chosen, as applications can prove.

Let us denote by $\tilde{C}$ the PET experimental concentration in a specific tissue, by $C_b^{\text{ROI}}$ the IF, and by $V_b$ the volume fraction occupied by blood. In principle, $V_b$ may depend on the position within the tissue but, with a good approximation, it can be considered constant since the blood perfusion, under physiological conditions, is homogeneous inside the same organ. 
Consider a compartmental model with $m$ compartments $1, \dots, m$ and an arbitrary number $p$ of exchange coefficients $k_1, \dots, k_p$. Then $\tilde{C}$ can be expressed as \begin{equation} \label{eqn:inv_pb1}
\tilde{C}(t) = (1-V_b)\sum_{i=1}^m C_i(t) + V_b C_b^{\text{ROI}}(t) \ \forall t \in \mathbb{R}_+ \ .
\end{equation}
Equivalently, we can write 
\begin{equation} \label{eqn:inv_pb2}
\tilde{C}(t) - V_b C_b^{\text{ROI}}(t) =  \boldsymbol{\alpha}^T \boldsymbol{C}(t) \ , 
\end{equation}
where 
\begin{equation*} 
\boldsymbol{\alpha} = \pmatrix{ 1-V_b \cr 1-V_b \cr \vdots \cr 1-V_b} \ \mbox{and} \ \boldsymbol{C} = \pmatrix{ C_1 \cr C_2 \cr \vdots \cr C_m } \ ,
\end{equation*}
and where we assembled the available experimental data at the left-hand side. 
In particular, the inverse problem equations for the two compartmental models considered in this paper are the following.
\begin{itemize}
\item Two--compartment catenary system: 
\begin{equation} \label{eqn:2C_inv_pb}
\hspace{-0.85cm}\tilde{C}(t) - V_b C_b^{\text{ROI}}(t) =  \boldsymbol{\alpha}^T \boldsymbol{C}(t) \ \mbox{with} \ \boldsymbol{\alpha} = \pmatrix{1-V_b \cr 1-V_b} \ \mbox{and} \ \boldsymbol{C} = \pmatrix{ C_f \cr C_m},
\end{equation}
where $\boldsymbol{C}$ is given by equation (\ref{eqn:2C_sol}).
\item Three--compartment non--catenary system: 
\begin{equation} \label{eqn:kidney_inv_pb}
\hspace{-0.85cm}\tilde{C}(t) - V_b C_a^{\text{ROI}}(t) =  \boldsymbol{\alpha}^T \boldsymbol{C}(t) \ \mbox{with} \ \boldsymbol{\alpha} = \pmatrix{1-V_b \cr 1-V_b \cr 1-V_b} \ \mbox{and} \ \boldsymbol{C} = \pmatrix{ C_f \cr C_m \cr C_t},
\end{equation}
where $\boldsymbol{C}$ is given by equation (\ref{eqn:kidney_sol2}).
\end{itemize}

We re-write equation (\ref{eqn:inv_pb2}) for the unknown vector parameter $\boldsymbol{k} = (k_1, \dots, k_p)$ in the form
\begin{equation} \label{eqn:inv_pb_zero}
\boldsymbol{\alpha}^T \boldsymbol{C}(t) + V_b C_b^{\text{ROI}}(t) - \tilde{C}(t) := \mathcal{F}_t(\boldsymbol{k}) = 0 \ ,
\end{equation} 
where $\mathcal{F}_t \colon \mathbb{R}_+^p \to C^1(\mathbb{R}_+,\mathbb{R})$ is a non--linear operator parameterized by the time variable $t \in \mathbb{R}_+$, whose dependence on the unknown vector $\boldsymbol{k}$ is made explicit. 
 The non--linear zero--finding problem of equation (\ref{eqn:inv_pb_zero}) is solved by means of the Gauss--Newton method, that reads
\begin{equation} \label{eq:GN}
\bigg[ \frac{d \mathcal{F}_t}{d \boldsymbol{k}} (\boldsymbol{k}^{(0)};\boldsymbol{\delta}^{(0)}) \bigg](t) = - \mathcal{F}_t(\boldsymbol{k}^{(0)}) \ ,
\end{equation}
with unknown step--size $\boldsymbol{\delta}^{(0)} \in \mathbb{R}^p$, initial guess $\boldsymbol{k}^{(0)} \in \mathbb{R}_+^p$ and for $t \in \mathbb{R}_+$. The operator $\mathcal{F}_t$ is differentiable and even analytic, therefore it is possible to compute analytically its Fr\'echet derivative, which is the bounded and linear operator 
\begin{eqnarray*} \label{eq:dFt_dk_def}
\frac{d \mathcal{F}_t}{d \boldsymbol{k}} \colon & \mathbb{R}^p \to C^1(\mathbb{R}_+,\mathbb{R}) \\
& \ \boldsymbol{\delta} \mapsto \bigg[t \mapsto \nabla_{\boldsymbol{k}} \mathcal{F}_t(\boldsymbol{k}) \cdot \boldsymbol{\delta} \bigg] \ .
\end{eqnarray*}  
In real applications, only noisy versions of $\tilde{C}(t)$ and $C_b^{\text{ROI}}$ for a finite number of sampling time points $t_1, \dots, t_T \in \mathbb{R}_+$ are available. Therefore, equation (\ref{eq:GN}) becomes the discretized linear system
\begin{equation} \label{eq:GN_discr}
\boldsymbol{F}^{(0)} \boldsymbol{\delta}^{(0)} = \boldsymbol{Y}^{(0)} \ ,
\end{equation} 
where
\begin{equation*} 
\boldsymbol{F}^{(0)} = \pmatrix{ \nabla_{\boldsymbol{k}} \mathcal{F}_{t_1}(\boldsymbol{k}^{(0)})^T \cr \vdots \cr \nabla_{\boldsymbol{k}} \mathcal{F}_{t_T}(\boldsymbol{k}^{(0)})^T } \ , \
\boldsymbol{Y}^{(0)} = \pmatrix{\tilde{C}(t_1) - V_b C_b^{\text{ROI}}(t_1) - \boldsymbol{\alpha}^T \boldsymbol{C}(t_1) \cr \vdots \cr \tilde{C}(t_T) - V_b C_b^{\text{ROI}}(t_T) - \boldsymbol{\alpha}^T \boldsymbol{C}(t_T)} \ .
\end{equation*} 
The system (\ref{eq:GN_discr}), with the step--size vector $\boldsymbol{\delta}^{(0)}$ as unknown, constitutes a classic linear ill--posed inverse problem, since the solution may not exist, may not be unique, and may not be stable. Regularization is needed in order to find a unique stable solution of (\ref{eq:GN_discr}). We consider a Tikhonov--type regularization \cite{Tikhonov}, with the Tikhonov penalty on the step--size vector, which leads to the regularized system  
\begin{equation} \label{eq:GN_discr_reg}
(r^{(0)} \boldsymbol{I}_{[p]} + {\boldsymbol{F}^{(0)}}^{T} \boldsymbol{F}^{(0)}) \boldsymbol{\delta}^{(0)} = {\boldsymbol{F}^{(0)}}^{T}\boldsymbol{Y}^{(0)} \ ,
\end{equation} 
where $r^{(0)}$ is the regularization parameter, which is allowed to change at every iteration, and $\boldsymbol{I}_{[p]}$ is the $p \times p$ identity matrix. The step--size $\boldsymbol{\delta}^{(0)}$ is the least--square solution of (\ref{eq:GN_discr_reg}). The optimization algorithm performs an iterative scheme  which increases the values of the exchange coefficients by letting $\boldsymbol{k}^{(1)} = \boldsymbol{k}^{(0)} + \boldsymbol{\delta}^{(0)}$, and iterates the process until a selected stopping criterion is satisfied. Notice that solving the problem by means of Tikhonov regularization is equivalent to asking for the solution of the original system (\ref{eq:GN_discr}) with a small norm, i.e. limiting the step--size length. This property allows to avoid divergence of the iterative algorithm. The role of the regularization parameter is crucial, since in general it supervises the importance of the regularization term, and in particular it defines the direction along which look for the solution. For example, for large values of the regularization parameter the step--size is taken approximately in the direction of the gradient.

The implementation details of the pixel--wise regularized Gauss--Newton algorithm are given in the forthcoming description of the parametric imaging method.

\subsection{Uniqueness theorems} \label{subsec:id}

The identifiability analysis \cite{Miao,Yates} of compartmental models is a necessary step in the solution process, since the uniqueness of the reconstructed parameters is an essential property which makes the model effective in the description of the physiological processes under investigation. 

As shown in \cite{Delbary_id}, the two--compartment catenary systems describing the FDG cellular metabolism is always identifiable, i.e. the following theorem holds.

\begin{theorem}
The inverse problem for the two--compartment catenary system of equation (\ref{eqn:2C_inv_pb}) 
has a unique solution $\boldsymbol{k} = (k_{fb}, k_{bf}, k_{mf}, k_{fm})^T \in \mathbb{R}_+^{4*}$  
determined from the PET experimental measurements $C_b^{\text{ROI}}$ and $\tilde{C}$.
\end{theorem}

We show here that the three--compartment non--catenary model for the renal system is also identifiable. The proof follows the idea presented in \cite{Delbary_id} of using the Laplace transform of the compartmental system of ODEs as the key mathematical tool, but here we consider a non--catenary model in a pixel--wise framework. The main novelty is in the key role of the modeling assumptions that $C_t(t)$, $k_{tm}$, $k_{ut}$ are fixed, which were obtained in \ref{subsubsec:ktm_kut} by extension to single pixels of physiological properties holding for the compartmental model of the entire renal system. These assumptions are essential for uniqueness. 

\begin{theorem}
 Let the concentration $C_t(t)$ and parameters $k_{tm}$ and $k_{ut}$ be fixed.
By assuming that the polynomials
\begin{equation*}
P(s) = k_{ma} (s+k_{af}+k_{mf}) + k_{fa} k_{mf} 
\end{equation*}
and
\begin{equation*}
Q(s) = k_{fa} (s+k_{fm}+k_{tm}) + k_{ma} k_{fm}
\end{equation*}
are both coprime with the polynomial
\begin{equation*}
D(s) = (s+k_{af}+k_{mf})(s+k_{fm}+k_{tm}) - k_{mf} k_{fm} \ ,
\end{equation*}
the inverse problem for the three--compartment non--catenary system of equation (\ref{eqn:kidney_inv_pb}) has a unique solution $\boldsymbol{k} =(k_{fa}, k_{ma}, k_{af}, k_{mf}, k_{fm}, k_{tm}, k_{ut})^T \in \mathbb{R}_+^{7*}$ determined from the PET experimental measurements $C_a^{\text{ROI}}$ and $\tilde{C}$. 
 
\begin{proof}
 To simplify notations we let
\begin{equation*}
k_{tm} = \alpha \ , \ k_{ut} = \beta \ ,
\end{equation*}
\begin{equation*}
k_{fa} = a \ , \ k_{ma} = b \ , \ k_{af} = c \ , \ k_{mf} = d \ , \ k_{fm} = e \ .
\end{equation*}
Since $k_{tm}$ and $k_{ut}$, i.e. $\alpha$ and $\beta$, are assumed to be fixed, we have to prove uniqueness of the remaining five coefficients $a, b, c, d, e$.
In the new notations the system of ODEs (\ref{eqn:kidney_ode}) takes the form
\begin{equation} \label{eqn:ode_id} 
\cases{\dot{C}_f = -(c + d) C_f + e C_m + a C_a^{\text{ROI}} \\
\dot{C}_m = d C_f - (e + \alpha) C_m + b C_a^{\text{ROI}} \\
\dot{C}_t = \alpha C_m - \beta C_t} \ .
\end{equation}

Assuming that the concentrations are sufficiently regular, we take the Laplace transform of the differential equations (\ref{eqn:ode_id}):
\begin{equation} \label{eqn:ode_id_laplace1}
\cases{
(s + c + d) \mathcal{L}(C_f) - e \mathcal{L}(C_m) = a \mathcal{L}(C_a^{\text{ROI}}) \\ 
-d \mathcal{L}(C_f) + (s + e + \alpha) \mathcal{L}(C_m) = b \mathcal{L}(C_a^{\text{ROI}}) \\
(s + \beta) \mathcal{L}(C_t) - \alpha \mathcal{L}(C_m) = 0} \ ,
\end{equation}
where $\mathcal{L}(f)$ denotes the Laplace transform of the function $f$.
Solving  (\ref{eqn:ode_id_laplace1}) with respect to $\mathcal{L}(C_f)$, $\mathcal{L}(C_m)$, and $\mathcal{L}(C_t)$, we get
\numparts
\begin{eqnarray}
\label{eqn:Cf_id_laplace}
\mathcal{L}(C_f) = \frac{a (s + e + \alpha) + b e}{D} \mathcal{L}(C_a^{\text{ROI}}) \ , \\
\label{eqn:Cm_id_laplace}
\mathcal{L}(C_m) = \frac{b (s + c + d) + a d}{D} \mathcal{L}(C_a^{\text{ROI}}) \ , \\
\label{eqn:Ct_id_laplace}
\mathcal{L}(C_t) = \frac{\alpha}{s + \beta} \mathcal{L}(C_m) \ , 
\end{eqnarray}
\endnumparts
where $D = (s + c + d) (s + e + \alpha) - d e$.
Moreover, by comparing (\ref{eqn:Cm_id_laplace}) and (\ref{eqn:Ct_id_laplace}), we get
\begin{equation} \label{eqn:Ct_id_laplace_bis}
\mathcal{L}(C_t) = \frac{\alpha}{s + \beta} \frac{b (s + c + d) + a d}{D} \mathcal{L}(C_a^{\text{ROI}}) \ .
\end{equation}
Then, we take the Laplace transform of equation (\ref{eqn:kidney_inv_pb}):
\begin{equation} \label{eqn:L_kidney_inv_pb}
\frac{\mathcal{L}(\tilde{C}) - V_b \mathcal{L}(C_a^{\text{ROI}}) }{1-V_b}  = \mathcal{L}(C_f) +\mathcal{L}(C_m) + \mathcal{L}(C_t) \ ,
\end{equation}
where the left-hand side is a known function of $s$, independent of the constants $a,b,c,d,e$.

Now, suppose $(a', b', c', d', e')$ is an alternative choice of rate coefficients consistent with the data of the problem.  We denote by $f'$ the correspondent of any function $f(a,b,c,d,e)$, with argument $(a,b,c,d,e)$ replaced by $(a',b',c',d',e')$.
The condition $C_t$ fixed implies that $C_t= C'_t$, whence it follows that
\begin{equation} \label{eqn:L_C_t}
\mathcal{L}(C_t) = \mathcal{L}(C'_t) \ .
\end{equation}
Equation (\ref{eqn:Ct_id_laplace}) implies
\begin{equation*} 
\mathcal{L}(C_m) = \frac{s + \beta}{\alpha} \mathcal{L}(C_t) \ , 
\end{equation*}
and, since $\alpha$ and $\beta$ are fixed, and equation (\ref{eqn:L_C_t}) holds, it follows also that 
\begin{equation} \label{eqn:L_C_m}
\mathcal{L}(C_m) = \mathcal{L}(C'_m) \ . 
\end{equation}
From equation (\ref{eqn:L_kidney_inv_pb}), it is found that
\begin{equation} \label{eqn:id_kidney_inv_pb}
\mathcal{L}(C_f) +\mathcal{L}(C_m) + \mathcal{L}(C_t)= \mathcal{L}(C'_f) +\mathcal{L}(C'_m) + \mathcal{L}(C'_t) \ . 
\end{equation}
Therefore, because of (\ref{eqn:L_C_t}) and (\ref{eqn:L_C_m}), equation (\ref{eqn:id_kidney_inv_pb}) reduces to 
\begin{equation} \label{eqn:L_C_f}
\mathcal{L}(C_f) =\mathcal{L}(C'_f) \ . 
\end{equation}
Now, substitution of (\ref{eqn:Ct_id_laplace_bis}) into (\ref{eqn:L_C_t}) leads to 
\begin{equation} \label{eqn:id1}
\frac{b (s + c + d) + a d}{D} = \frac{b' (s + c' + d') + a' d'}{D'} \ ,
\end{equation}
and substitution of (\ref{eqn:Cf_id_laplace}) into (\ref{eqn:L_C_f}) gives
\begin{equation} \label{eqn:id2}
\frac{a (s + e + \alpha) + b e}{D} = \frac{a' (s + e' + \alpha) + b' e'}{D'} \ .
\end{equation}

Assuming that the two rational fractions (\ref{eqn:id1}) and (\ref{eqn:id2}) are irreducible, i.e. the polynomials $P(s) = b (s + c + d) + a d$ and $Q(s) = a (s + e + \alpha) + b e$ are both coprime with $D(s)$, we obtain the links between the two sets of parameters, i.e. the system
\begin{equation*}
\cases{
b = b' \\
b (c + d) + a d = b' (c' + d') + a' d' \\
c + d + e + \alpha = c' + d' + e' + \alpha \\
c e + (c + d) \alpha = c' e' + (c' + d')  \alpha \\
a = a' \\
a (e + \alpha) + b e = a' (e' + \alpha) + b' e' 
}
\end{equation*} 
which holds if and only if
\begin{equation*}
a = a' \ , \ b = b' \ , \ c = c' \ , \ d = d' \ , \ e = e' \ .
\end{equation*}
\end{proof}
\end{theorem}

\section{The imaging method} \label{sec:imaging}

In this section we present a parametric imaging method, which relies upon the application of image pre--processing algorithms and of a rather general optimization scheme based on the regularized Gauss--Newton method, presented in Section~\ref{sec:inverse}, able to solve the compartmental inverse problem pixel--wise. We remark that our method is potentially applicable to any generic compartmental model, provided an \emph{ad hoc} identifiability study and taking into account the compartment--dependent increase of the computational cost.

We start from the set of $N$ reconstructed dynamic FDG--PET images:
\begin{equation}
(\boldsymbol{f}_1^{(t)}, \boldsymbol{f}_2^{(t)}, \dots, \boldsymbol{f}_N^{(t)}) \mbox{ for } t = 1, \dots, T \ ,
\end{equation}
where $\boldsymbol{f}_n^{(t)}$ is the $n-$th PET image at $t-$th time point of tracer concentration $\tilde{C}$, i.e.
\begin{equation}
\boldsymbol{f}_n^{(t)}(i,j) = \tilde{C}_{(i,j)}(t) \mbox{ for } i = 1, \dots , I , \ j = 1, \dots, J \ ,
\end{equation}
and $I, \ J$ are the image dimensions.

We select the tissue of interest and the compartmental model reliable for its functional description. 
For each dynamic PET image $(\boldsymbol{f}_{\bar{n}}^{(1)}, \dots, \boldsymbol{f}_{\bar{n}}^{(T)})$, $\bar{n} \in \{ 1, \dots , N \}$, i.e. a PET slice, our imaging method follows the steps described below.

\begin{itemize}

\item[Step 1.] \emph{Gaussian smoothing.} In order to reduce the noise due to data acquisition, we apply a  truncated Gaussian smoothing filter through the convolution operation
\begin{equation}
\tilde{\boldsymbol{f}}_{\bar{n}}^{(t)} = \boldsymbol{f}_{\bar{n}}^{(t)} \ast \boldsymbol{G}_{0,\sigma} \ \forall t = 1, \dots, T \ , 
\end{equation}
where 
\begin{eqnarray}
\boldsymbol{G}_{0,\sigma}(i,j) = \frac{1}{2 \pi \sigma^2} e^{-{\frac{\boldsymbol{x}(i,j)^2+\boldsymbol{y}(i,j)^2}{2 \sigma^2}}} \mbox{ and } \\
\boldsymbol{x}(i,j) , \boldsymbol{y}(i,j) \in \{-\frac{L-1}{2}, \dots, \frac{L-1}{2}\} \times \{-\frac{L-1}{2}, \dots, \frac{L-1}{2}\} \ ,
\end{eqnarray}
and $L$ is the (odd) dimension of the window.
In all our applications, we use a Gaussian convolution matrix $\boldsymbol{G}_{0,\sigma}$ with zero mean, standard deviation $\sigma = 1$ and dimension $L=3$.  

\item[Step 2.] \emph{Image segmentation.} We model the outer region of the organ of interest with a standard two--compartment model; the ROI $\boldsymbol{A}$ delimiting the organ of physiologic interest is described by the most reliable compartmental model (according to the organ physiology). To extract the ROI $\boldsymbol{A}$ we apply the following image segmentation method:
\begin{enumerate}
\item compute the PET image averaged in time: $\displaystyle \tilde{\boldsymbol{f}}_{\bar{n}} = \frac{1}{T} \sum_{t=1}^T \tilde{\boldsymbol{f}}_{\bar{n}}^{(t)}$;
\item consider the pixel with maximum intensity: $(\bar{i},\bar{j}) = \max_{i,j} \tilde{\boldsymbol{f}}_{\bar{n}}(i,j)$;
\item approximate the profile of the $\bar{i}-$th matrix row with a family of one--dimensional Gaussian functions $G_{\bar{j},\sigma}$ of mean $\bar{j}$, variable variance $\sigma$,  and s.t. $G_{\bar{j},\sigma}(\bar{j})= \tilde{\boldsymbol{f}}_{\bar{n}}(\bar{i},\bar{j})$, by means of a curve fitting process. This consists in computing
\begin{equation}
\bar{\sigma}_{\bar{i}} = \underset{\sigma}{\mbox{arg}\min} || \tilde{\boldsymbol{f}}_{\bar{n}}(\bar{i},j) - G_{\bar{j},\sigma}(j) ||_{2} \ ;
\end{equation}
\item determine the activity's lower bound in the ROI as the value $\bar{c}$ at which the two curves $\tilde{\boldsymbol{f}}_{\bar{n}}(\bar{i},j)$ and  $G_{\bar{j},\bar{\sigma}_{\bar{i}}}(j)$ separate from each other. Formally, this consists in evaluating
\begin{equation}
\cases{
j^\ast = \underset{j \in (\bar{j}-\gamma,\bar{j}+\gamma)}{\mbox{arg}\max} | \tilde{\boldsymbol{f}}_{\bar{n}}(\bar{i},j) - G_{\bar{j},\bar{\sigma}_{\bar{i}}}(j) | \\
\bar{c} = \tilde{\boldsymbol{f}}_{\bar{n}}(\bar{i},j^\ast) 
} \ ,
\end{equation}
in a suitably chosen neighbourhood of $\bar{j}$ (i.e. for a suitable choice of $\gamma>0$). The ROI encompassing the organ is thus defined as
\begin{equation}
\boldsymbol{A}(i,j) = 
\cases{
0 \mbox{ if } \tilde{\boldsymbol{f}}_{\bar{n}}(i,j) < \bar{c} \\
1 \mbox{ if } \tilde{\boldsymbol{f}}_{\bar{n}}(i,j) \ge \bar{c} 
} \ .
\end{equation}
\end{enumerate}

\item[Step 3.] \emph{Parameter estimation.} We apply the regularized Gauss--Newton algorithm pixel-by-pixel, considering each pixel with its specific compartmental model. In general, for a compartmental model with $p$ arbitrary kinetic parameters, for a image pixel $(i,j) \in \{1, \dots , I\} \times \{1, \dots , J\}$, the reconstruction iterative algorithm reads as follows: 
\begin{enumerate}
\item check whether the radioactivity is significant: fix a constant value $\tau>0$ (e.g. $\tau = 10^2$) discriminating between background noise and tissue activity,
\begin{itemize}
\item if $|| \tilde{C}_{(i,j)}(t) ||_2 \le \tau$, then assign $\boldsymbol{k}=0$ and stop; 
\item if $|| \tilde{C}_{(i,j)}(t) ||_2 > \tau$, then continue;
\end{itemize}
\item choose the initial guess: $\boldsymbol{k}_{(i,j)}^{(0)} \in \mathbb{R}_+^p$;
\item solve for $\boldsymbol{\delta}_{(i,j)}^{(0)} \in \mathbb{R}^p$ 
\begin{equation}
(r^{(0)} \boldsymbol{I}_{[p]} + {\boldsymbol{F}^{(0)}_{(i,j)}}^{T} \boldsymbol{F}^{(0)}_{(i,j)} ) \boldsymbol{\delta}_{(i,j)}^{(0)} = {\boldsymbol{F}^{(0)}_{(i,j)}}^{T} \boldsymbol{Y}^{(0)}_{(i,j)} \ ,
\end{equation}
where $\boldsymbol{F}^{(0)}_{(i,j)}$  and $\boldsymbol{Y}^{(0)}_{(i,j)}$ are defined as in (\ref{eq:GN_discr}), and $r^{(0)}$ is the regularization parameter  automatically selected by means of the Generalized Cross Validation (GCV) method \cite{GCV,OSullivan}. The advantages in using the GCV are mainly that it can be applied without any {\it a priori} information on the error on the data or on peculiar properties of the solution, and that it requires just the computation of the SVD of the matrix of the problem;
\item project onto zero the components of the step--size vector $\boldsymbol{\delta}_{(i,j)}^{(0)}$ that make negative the components of the parameter vector $\boldsymbol{k}_{(i,j)}^{(0)}$. This means defining the $p \times p$ projection matrix $\boldsymbol{P}_{(i,j)}^{(0)}$ s.t.
\begin{equation}
\boldsymbol{P}_{(i,j)}^{(0)}(q,r) = 
\cases{
0 \mbox{ if } q \ne r \\ 
0 \mbox{ if } q = r \mbox{ and } (\boldsymbol{k}_{(i,j)}^{(0)})_q+(\boldsymbol{\delta}_{(i,j)}^{(0)})_q < 0 \\
1 \mbox{ if } q = r \mbox{ and } (\boldsymbol{k}_{(i,j)}^{(0)})_q+(\boldsymbol{\delta}_{(i,j)}^{(0)})_q > 0 
} \ ;
\end{equation}
\item update $\boldsymbol{k}_{(i,j)}^{(0)}$ with the projected step--size:
\begin{equation}
\boldsymbol{k}_{(i,j)}^{(1)} = \boldsymbol{k}_{(i,j)}^{(0)} + \boldsymbol{P}_{(i,j)}^{(0)} \boldsymbol{\delta}_{(i,j)}^{(0)} \ , 
\end{equation}
and iterate.
\end{enumerate}
The iterative scheme is stopped when the relative error between the experimental dynamic concentration and the model-predicted one is less than an appropriate threshold, i.e. at a generic iteration $h$
\begin{equation}
\frac{|| \tilde{C}_{(i,j)}(t) - V_b C_b^{\text{ROI}}(t) - \boldsymbol{\alpha}^T \boldsymbol{C}_{(i,j)}(t;\boldsymbol{k}_{(i,j)}^{(h)}) ||_2}{|| \tilde{C}_{(i,j)}(t)||_2} \le \epsilon \ ,
\end{equation}
where $\epsilon$ depends on the noise level on data.

\item[Step 4.] \emph{Parametric images.} Once we obtain the set of exchange coefficients of the model for each image pixel, we build up the parametric images $\boldsymbol{K}_1, \dots, \boldsymbol{K}_p$:
\begin{equation}
\boldsymbol{K}_1(i,j) = \boldsymbol{k}_{(i,j)}(1), \ \dots \ , \boldsymbol{K}_p(i,j) = \boldsymbol{k}_{(i,j)}(p) \ .
\end{equation}
\end{itemize}

Regarding Step 1. and Step 2., other imaging processing methods can be used to smooth and segment PET images; the impact of other approaches on the accuracy of compartmental analysis is under investigation. However, simple gaussian smoothing, as in Step 1., and our \emph{ad hoc} image segmentation process, as in Step 2., provide good results regardless of the  limited resolution of PET images involved in the analysis.

\section{Numerical applications}
Our new parametric imaging method is here applied against both synthetic and real microPET data of murine models \cite{microPET1, microPET2}. 
In particular, we test the procedure with the two--compartment model in the simulation framework and with the three--compartment model, describing the renal system, in the real data framework.

\subsection{Simulation}

\begin{figure}[htb]
\centering
\subfigure[\label{fig:phantom}]
{\includegraphics[width=3.7cm]{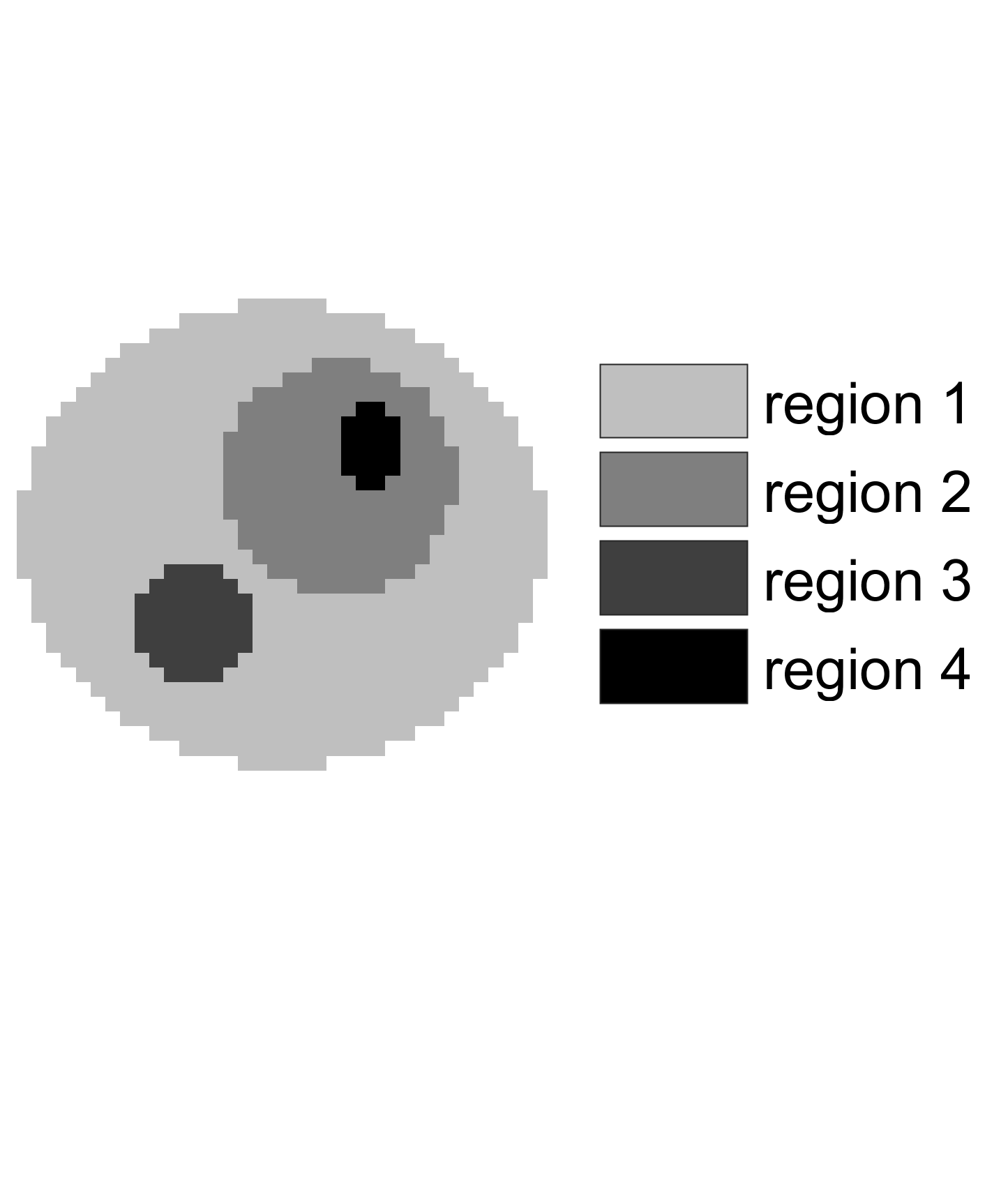}} \quad
\subfigure[\label{fig:IF_simulated}]
{\includegraphics[width=5cm]{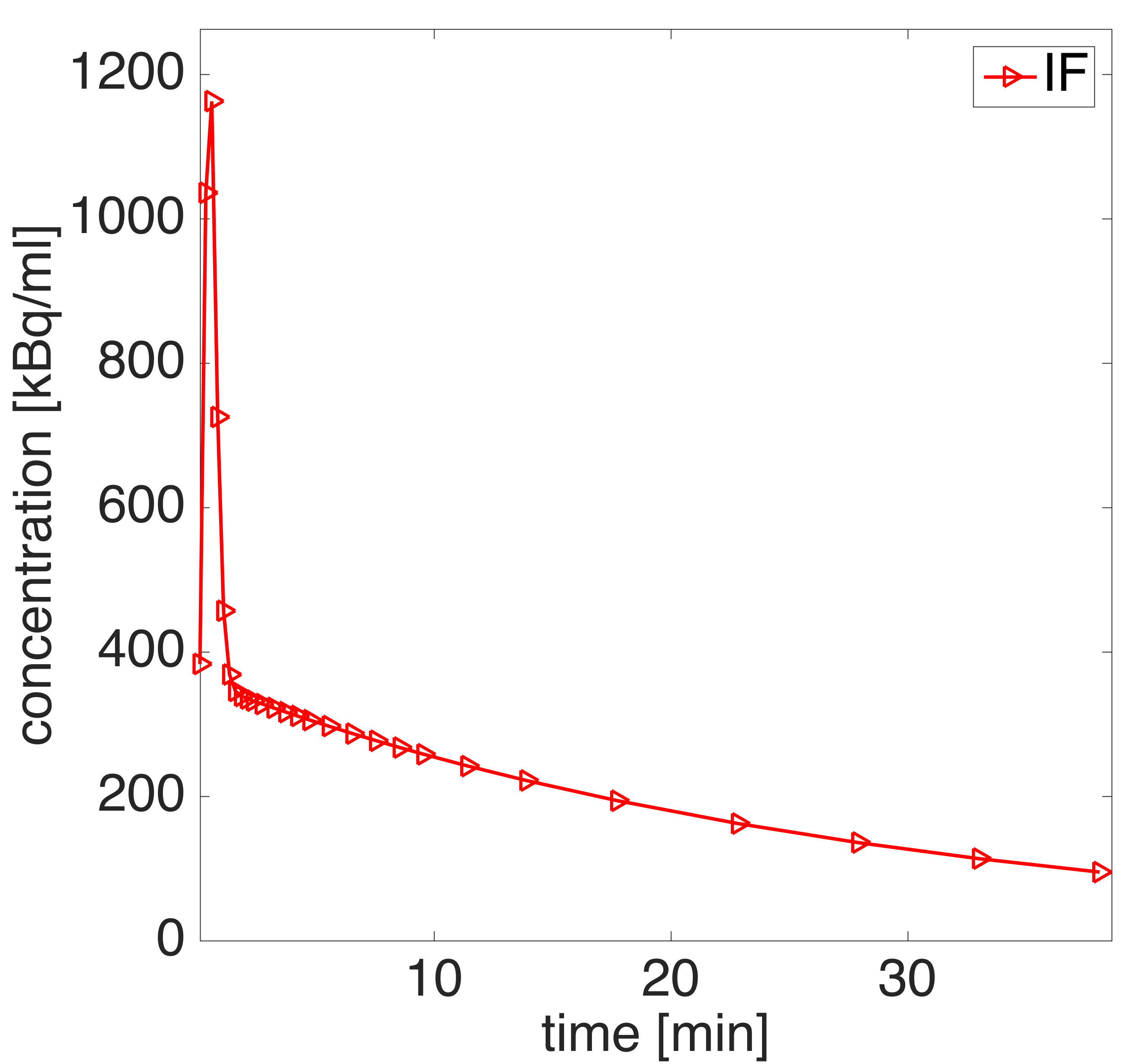}} \quad
\subfigure[\label{fig:TACs}]
{\includegraphics[width=5cm]{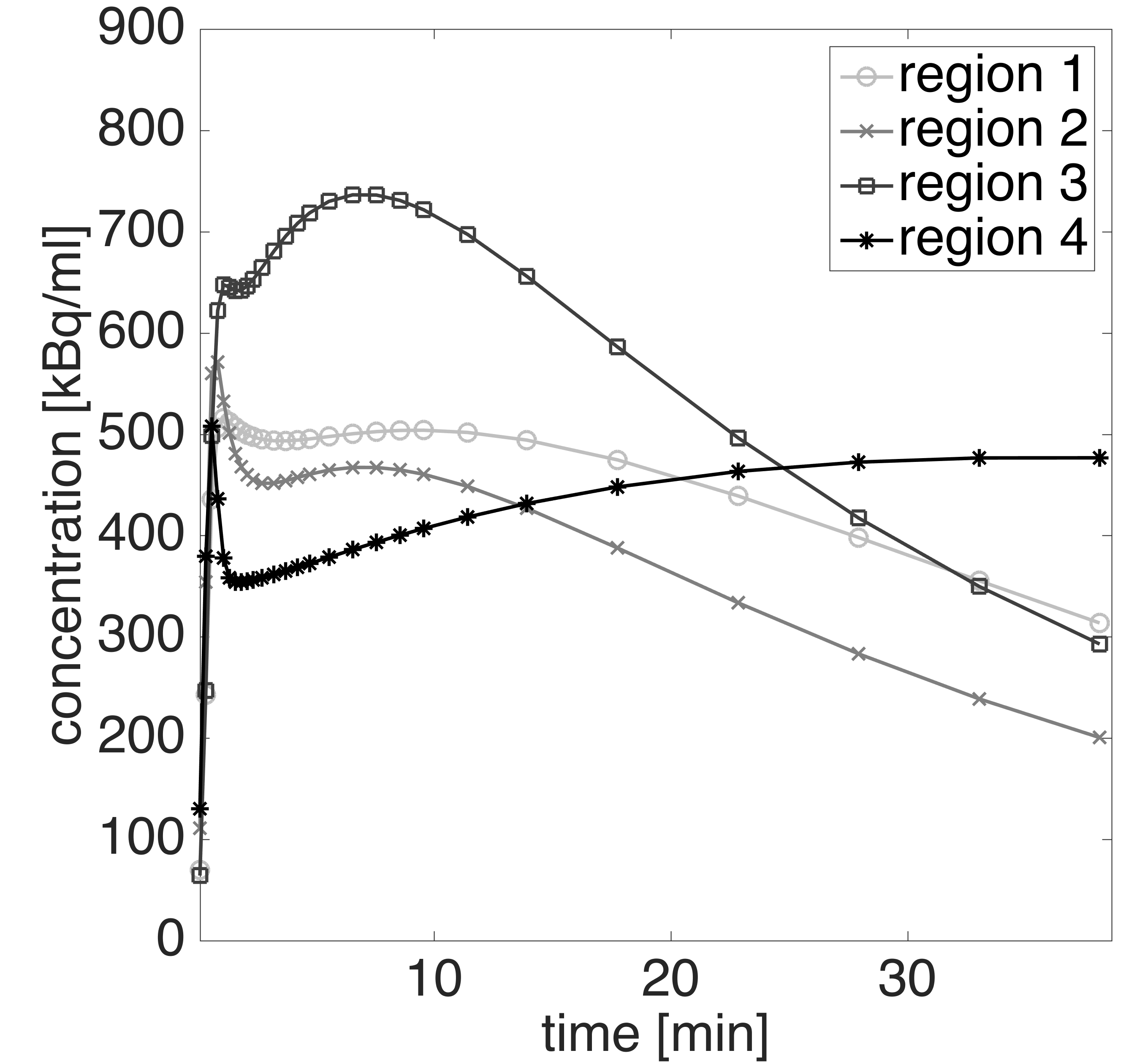}} \\
\subfigure[\label{fig:TACnoisy1}]
{\includegraphics[width=3.8cm]{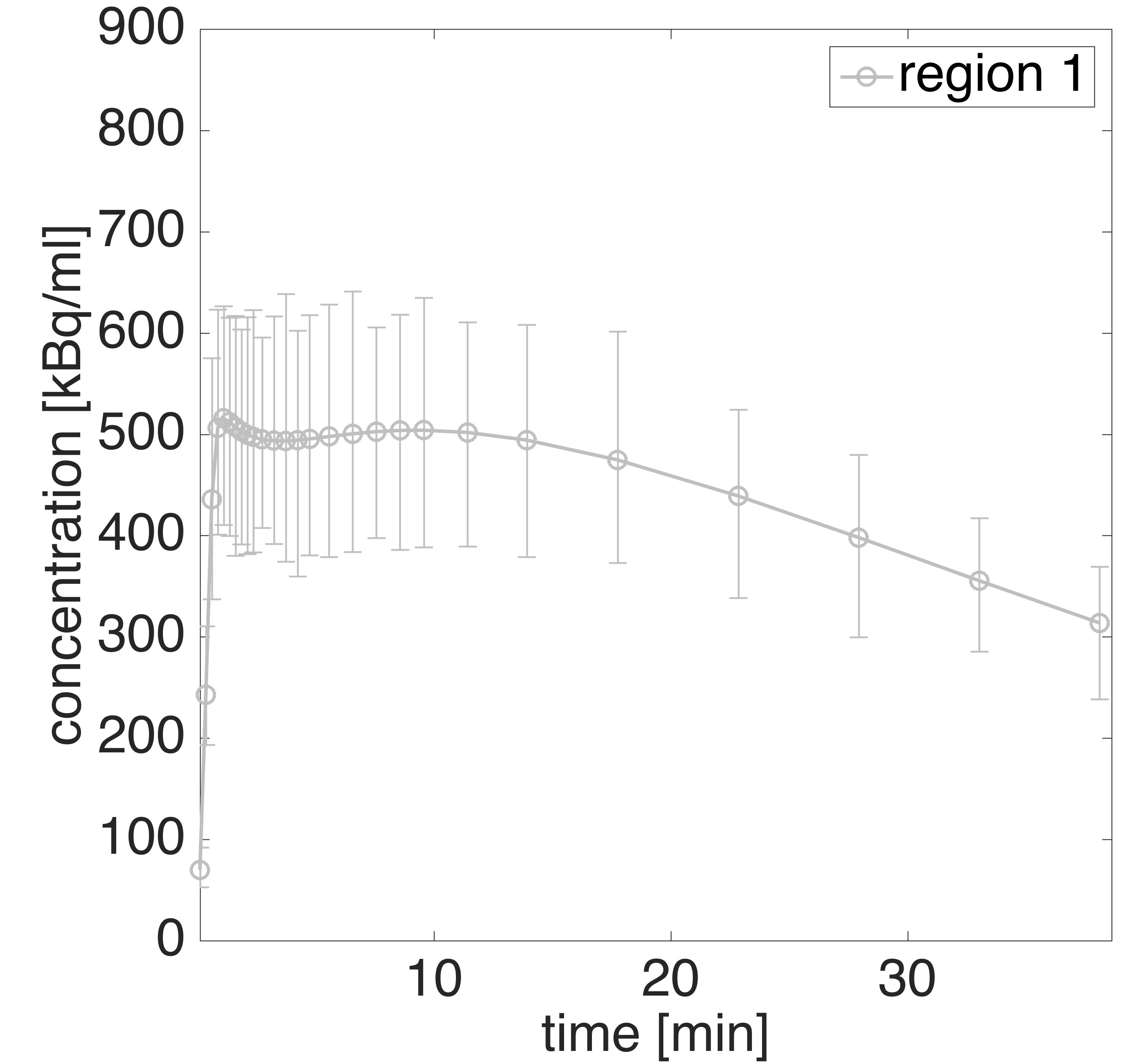}} 
\subfigure[\label{fig:TACnoisy2}]
{\includegraphics[width=3.8cm]{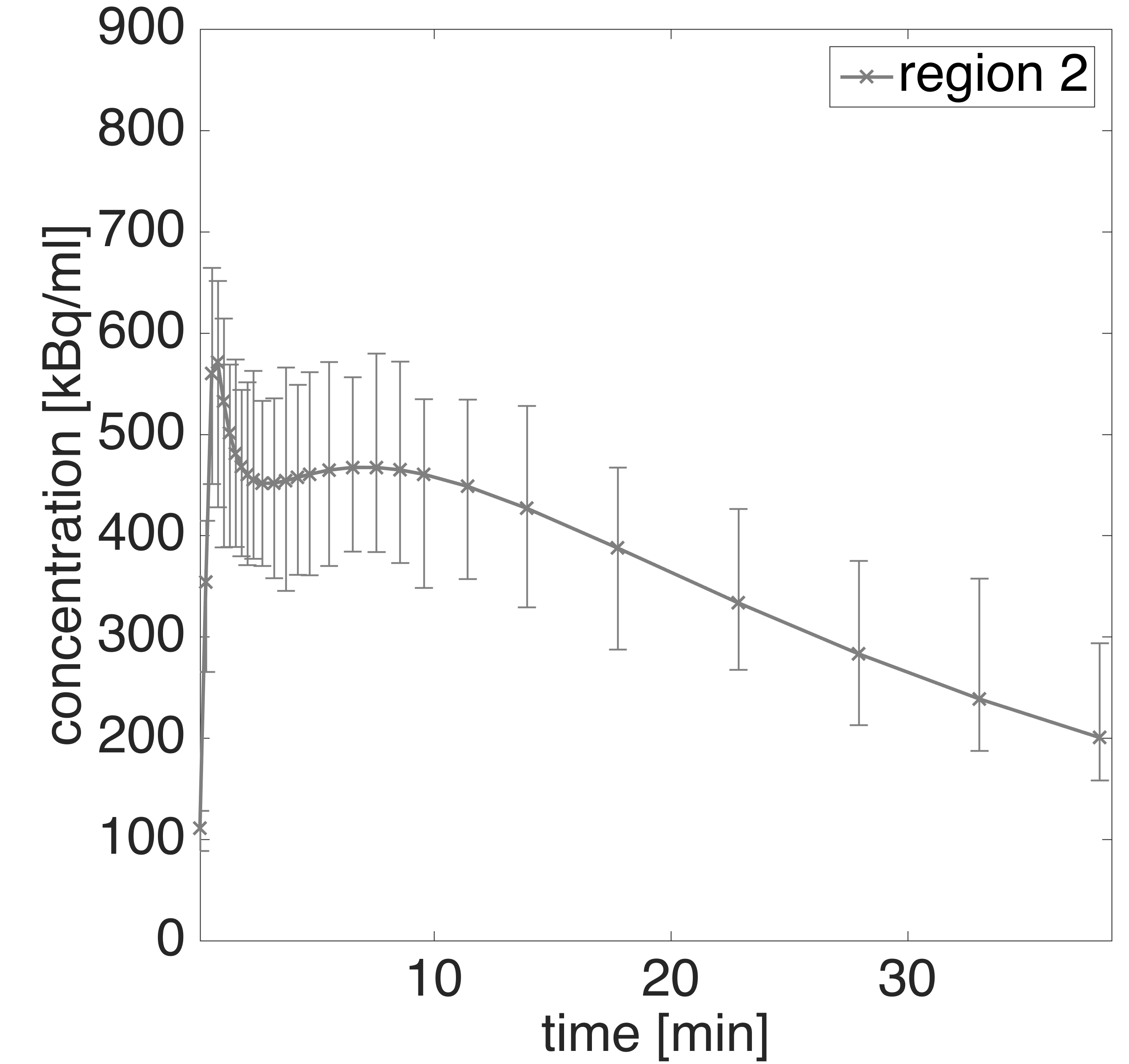}} 
\subfigure[\label{fig:TACnoisy3}]
{\includegraphics[width=3.8cm]{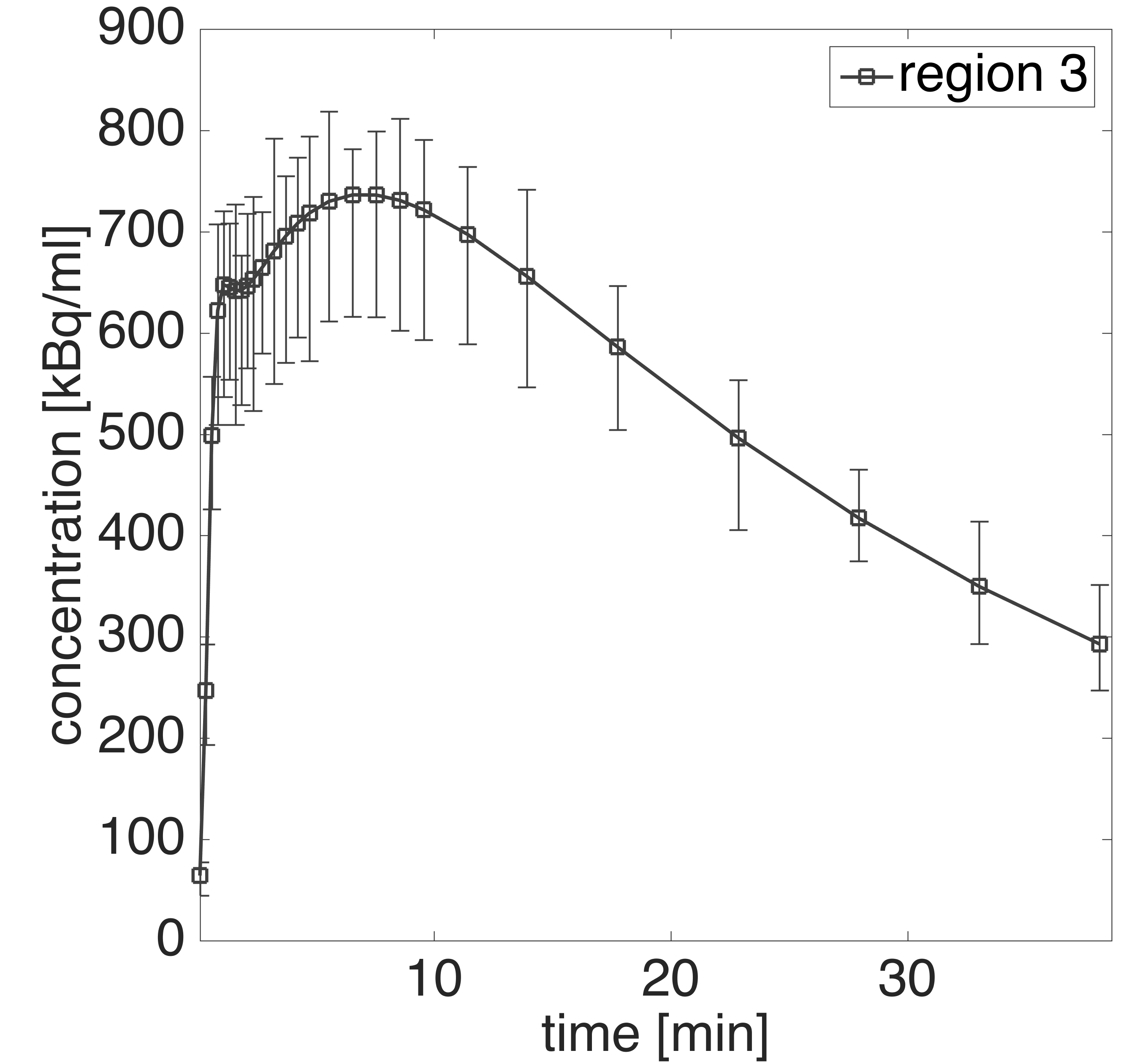}}
\subfigure[\label{fig:TACnoisy4}]
{\includegraphics[width=3.8cm]{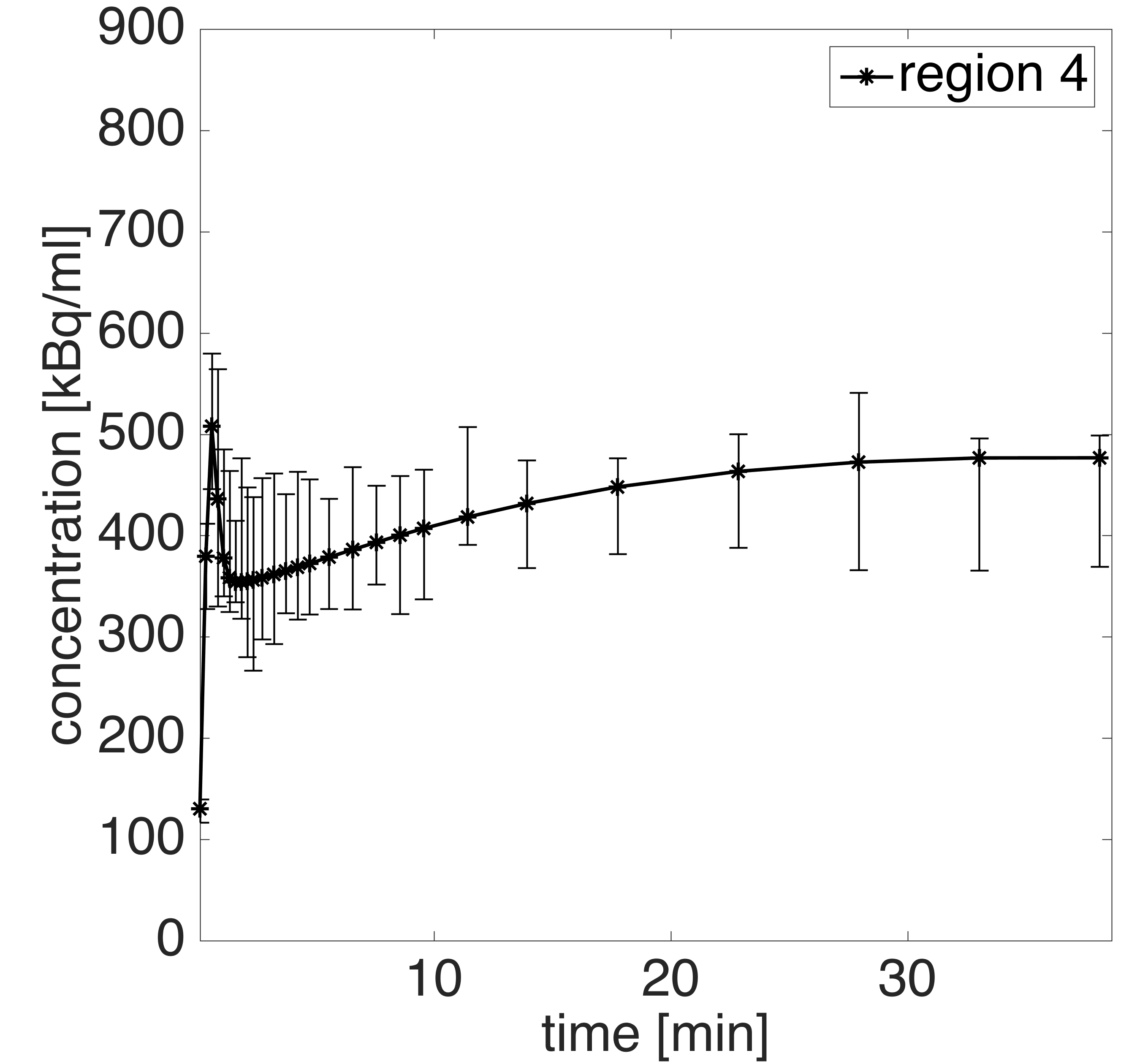}} 
\caption{The FDG--PET simulation setting. (a) Phantom composed by four homogeneous regions. (b) Simulated blood IF. (c) Characteristic noise--free TACs for the four regions. (d--g) Noisy TACs for each region: the error bars identify the variability on the concentrations introduced by the noise.}
\label{fig:simulation}
\end{figure}

\begin{figure}[htb]
\centering
\subfigure[$\boldsymbol{K}_{fb}$ \label{fig:K1}]
{\includegraphics[width=3.8cm]{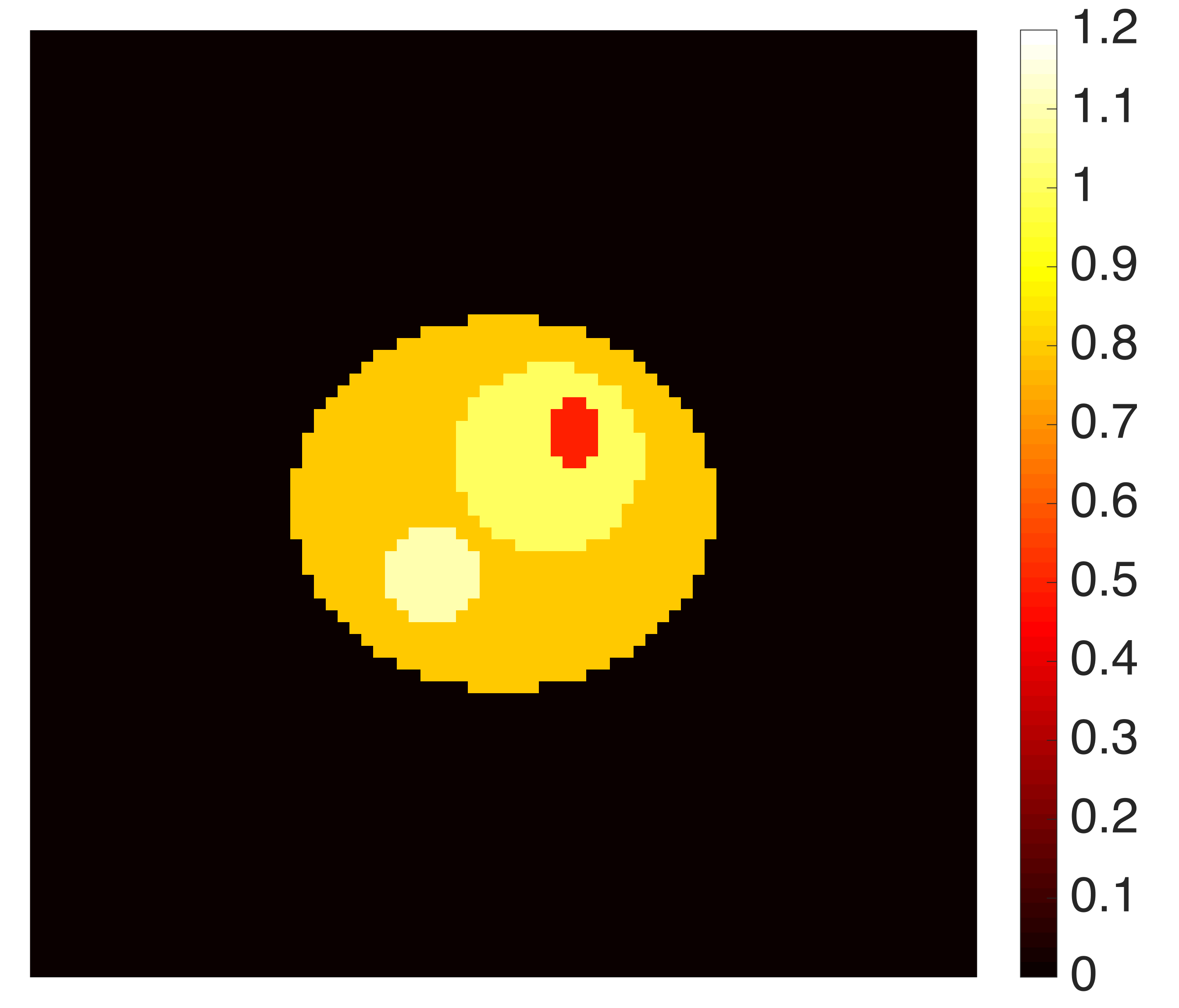}} 
\subfigure[$\boldsymbol{K}_{bf}$ \label{fig:K2}]
{\includegraphics[width=3.8cm]{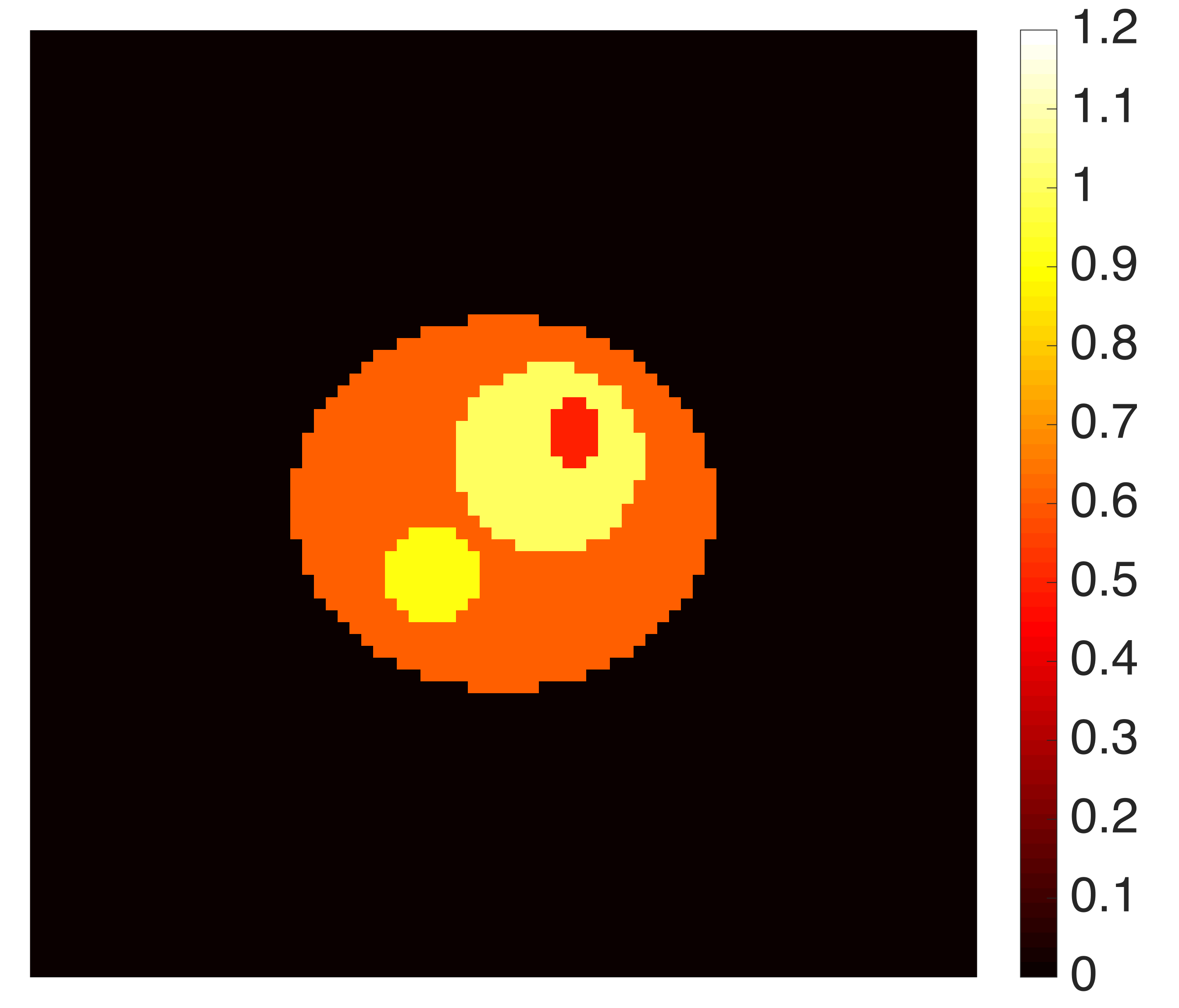}}
\subfigure[$\boldsymbol{K}_{mf}$ \label{fig:K3}]
{\includegraphics[width=3.8cm]{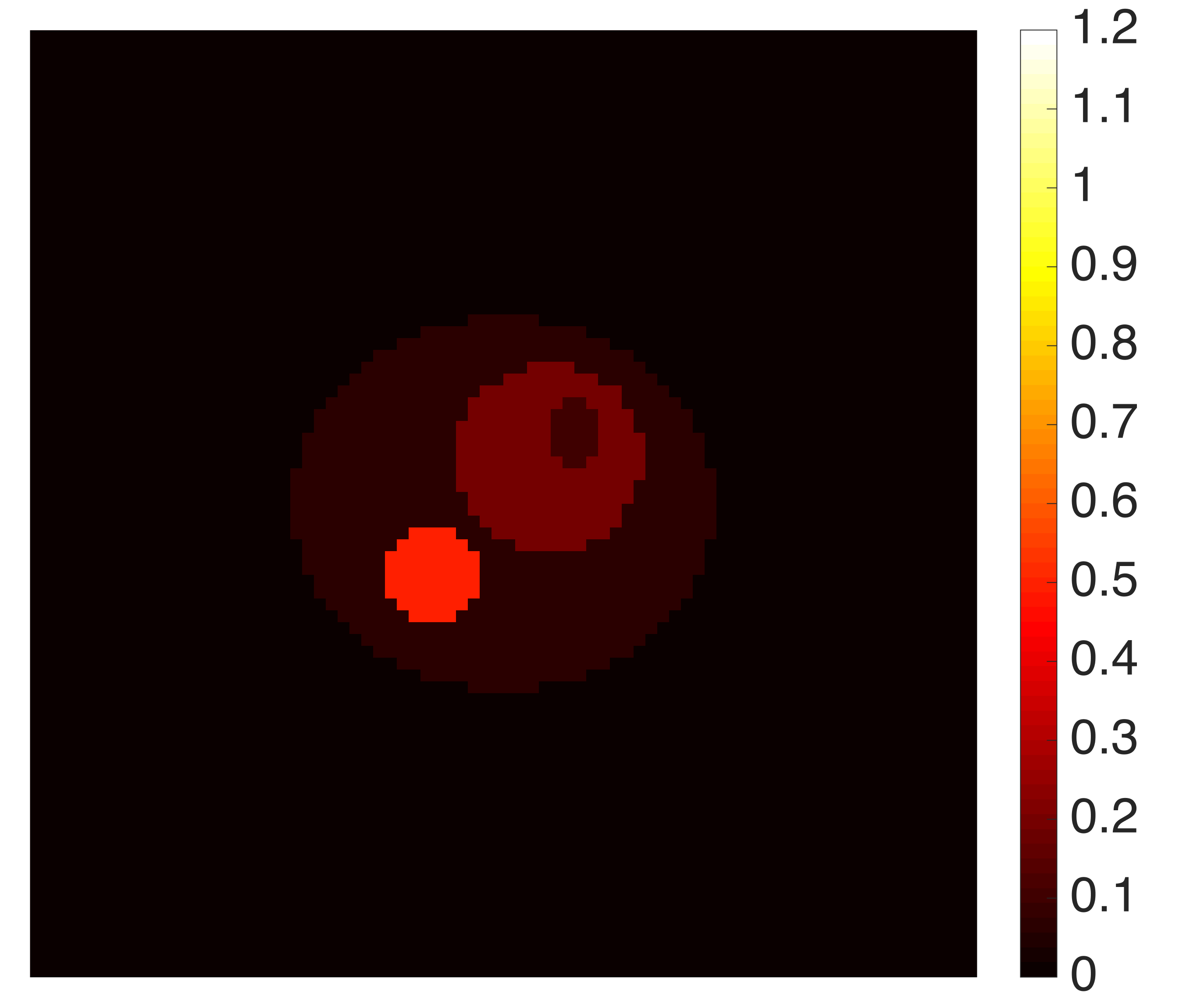}} 
\subfigure[$\boldsymbol{K}_{fm}$ \label{fig:K4}]
{\includegraphics[width=3.8cm]{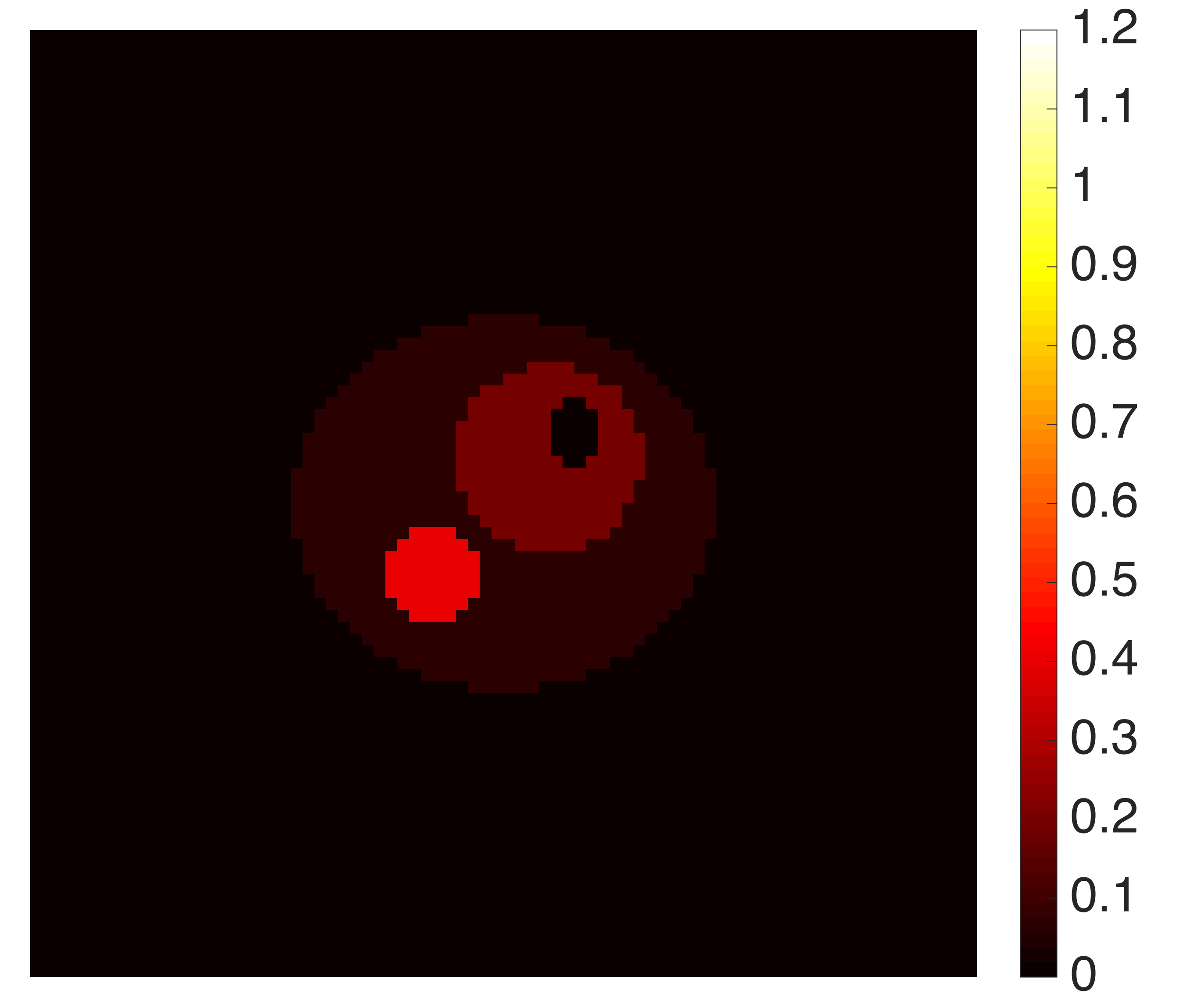}}
\caption{The ground truth parametric images $\boldsymbol{K}_{fb},\boldsymbol{K}_{bf},\boldsymbol{K}_{mf},\boldsymbol{K}_{fm}$ of the two--compartment system.}
\label{fig:parametric_images}
\end{figure}

Synthetic data are created mimicking a real FDG--microPET acquisition: first we chose a phantom (as in \figurename~\ref{fig:phantom}) encompassing four homogeneous regions; for every region, a set of realistic kinetic parameters of a two--compartmental problem was assigned as ground truth, and a realistic value for the blood volume fraction $V_b$ was selected
(\tablename~\ref{tab:kinetics}). 
Therefore, we obtained four synthetic parametric images  $\boldsymbol{K}_{fb},\boldsymbol{K}_{bf},\boldsymbol{K}_{mf},\boldsymbol{K}_{fm}$, each one characterized by a specific set of kinetic parameters. The ground truth parametric images are displayed in \figurename~\ref{fig:parametric_images}.

\begin{table}[htb] 
\caption{\label{tab:kinetics} Ground truth numerical values of the kinetic parameters $k_{fb},k_{bf},k_{mf},k_{fm}$ (min$^{-1}$), and of the blood volume fraction $V_b$, for each one of the four homogeneous regions.}
\lineup
\begin{indented}
\item[] \begin{tabular}{@{}cccccc}
\br
 & $k_{fb}$ & $k_{bf}$ & $k_{mf}$ & $k_{fm}$ & $V_b$ \\ 
\mr
region 1 & $0.8$ & $0.6$ & $0.07$ & $0.07$ & $0.1$ \\ 
region 2 & $1$ & $1$ & $0.2$ & $0.2$ & $0.2$ \\                             
region 3 & $1.1$ & $0.9$ & $0.5$ & $0.4$ & $0.05$ \\                             
region 4 & $0.5$ & $0.5$ & $0.1$ & $0.01$ & $0.3$ \\                           
\br
\end{tabular}
\end{indented}
\end{table} 

The dynamic PET images of tracer concentration were generated following the scheme below. For each pixel:
\begin{enumerate}
\item compute the values of the integral (\ref{eqn:2C_sol}) using the ground truth values of the parameters and a simulated blood IF (\figurename~\ref{fig:IF_simulated}), the latter obtained by fitting with a gamma variate function \cite{Golish} a set of real measurements acquired from a healthy mouse in a controlled experiment;
\item evaluate the synthetic concentration by means of equation (\ref{eqn:2C_inv_pb}), with the numerical value of $V_b$ related to the homogeneous region at which the pixel belongs;
\item sample on time interval $[t_1; t_T]$ of 27 time frames equivalent to the typical total acquisition time of the FDG experiments performed with the microPET scanner available at our lab (Albira, Carestream Health, Genova) \cite{Albira} and in agreement with usual time points of the experiments ($10 \times 15$s, $1 \times 22$s, $4 \times 30$s, $5 \times 60$s, $2 \times 150$s and $5 \times 300$s).
\end{enumerate} 
Once the noise--free dynamic PET images were obtained:
\begin{enumerate}
\item[(iv)] project the images into the sinogram space by means of the Radon transform, yielding projected noise--free sinogram data; 
\item[(v)] add mixed Poisson--Gaussian noise \cite{Luisier,Santarelli} to the projected data: apply Poisson noise to account for the stochastic nature of the photon counting process at the detectors, and then corrupt the Poisson model with additive Gaussian noise to account for the intrinsic thermal and electronic fluctuations of the acquisition device. Notice that the errors arising from instrumental and physical effects, such as attenuation, scattered events, decay and accidental coincidences, were not simulated; 
\item[(vi)] reconstruct the noisy dynamic PET images of tracer concentration by means of the Filtered Back Projection (FBP) applied on the noisy sinogram data. 
\end{enumerate} 
We created fifty independent identically--distributed noisy datasets. Characteristic noise--free TACs of the four regions are shown in \figurename~\ref{fig:TACs}, whereas the noisy TACs are reported in \figurename~\ref{fig:TACnoisy1}--\ref{fig:TACnoisy4}. The Poisson noise was applied by using the Matlab function \emph{poissrnd}, and the white Gaussian noise had a signal-to-noise ratio of 20 dB. 

For each dataset, we followed the reconstruction steps of Section~\ref{sec:imaging}, i.e. we applied the Gaussian smoothing filter ($\sigma = 1$, window $3 \times 3$) on the dynamic PET images and solved pixel--wise the compartmental inverse problem by means of the regularized Gauss--Newton algorithm. For each pixel, the starting point of our method was randomly chosen in the interval $(0 ,1)$ and the regularization parameter was optimized at each iteration through the GCV method.
We did not need to apply the image segmentation step because in this simulation we modeled the same two--compartment scheme for all the pixels.

\begin{figure}[htb]
\centering
\subfigure[$\boldsymbol{K}_{fb}$ \label{fig:K1}]
{\includegraphics[width=3.8cm]{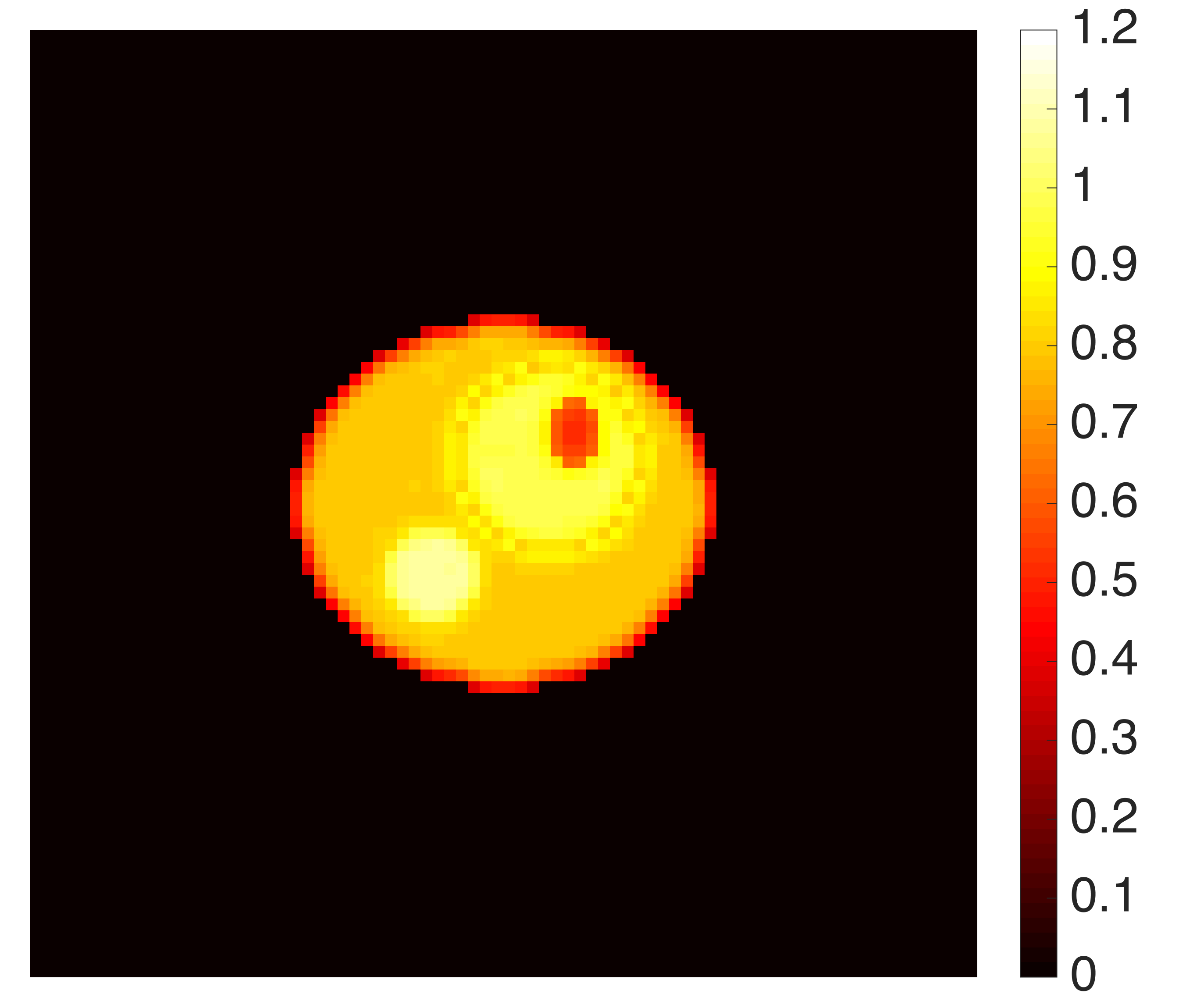}} 
\subfigure[$\boldsymbol{K}_{bf}$ \label{fig:K2}]
{\includegraphics[width=3.8cm]{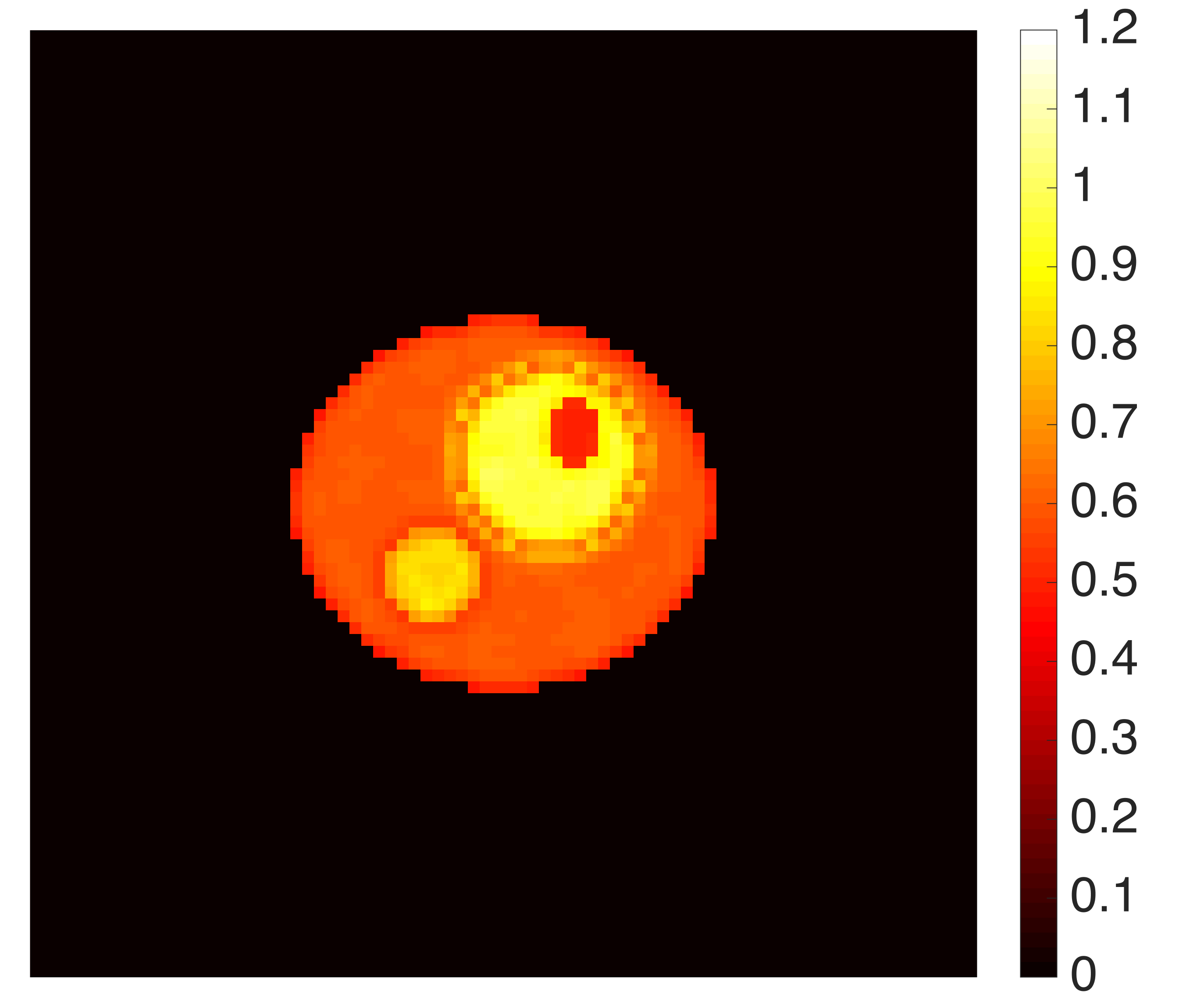}} 
\subfigure[$\boldsymbol{K}_{mf}$ \label{fig:K3}]
{\includegraphics[width=3.8cm]{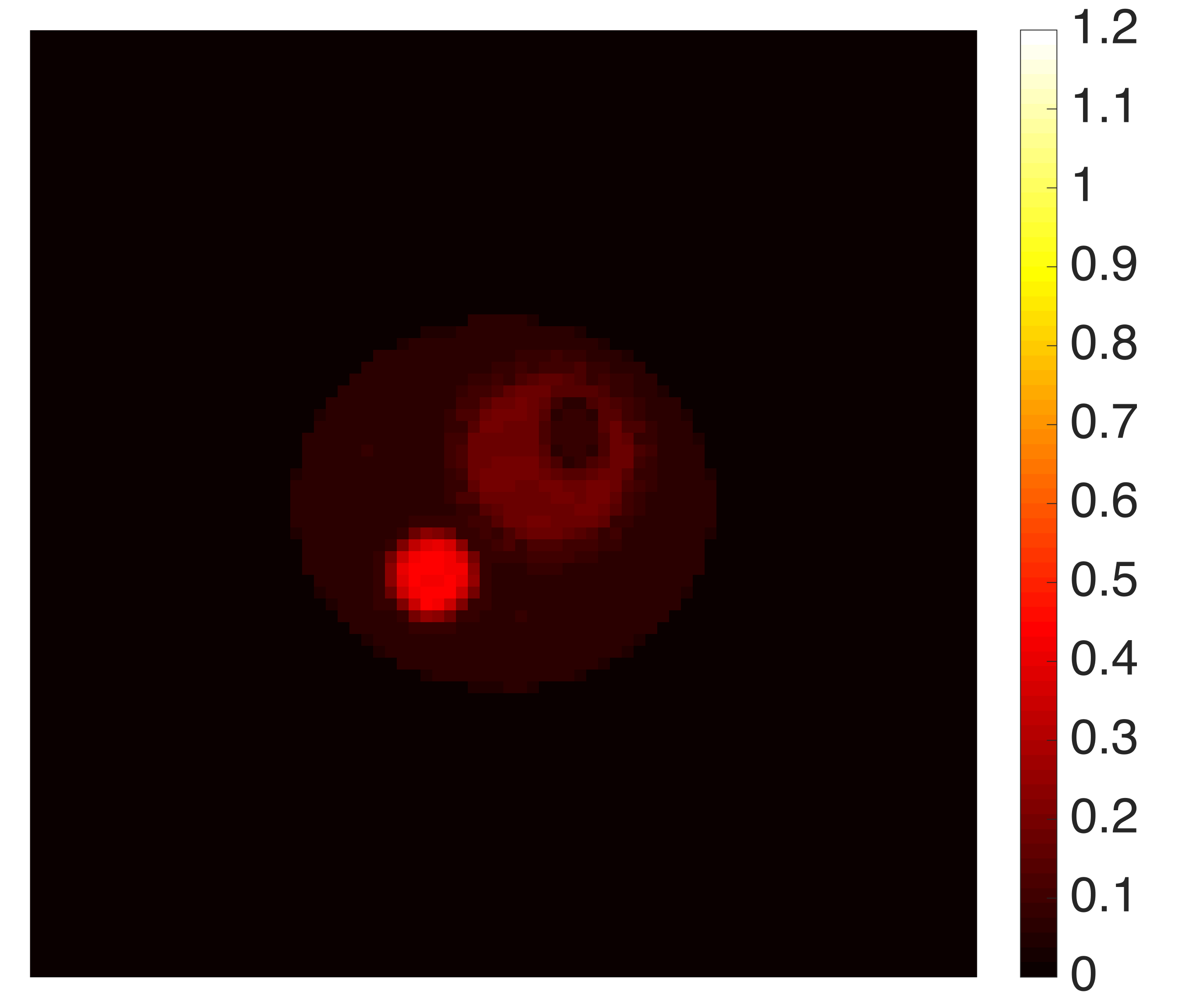}} 
\subfigure[$\boldsymbol{K}_{fm}$ \label{fig:K4}]
{\includegraphics[width=3.8cm]{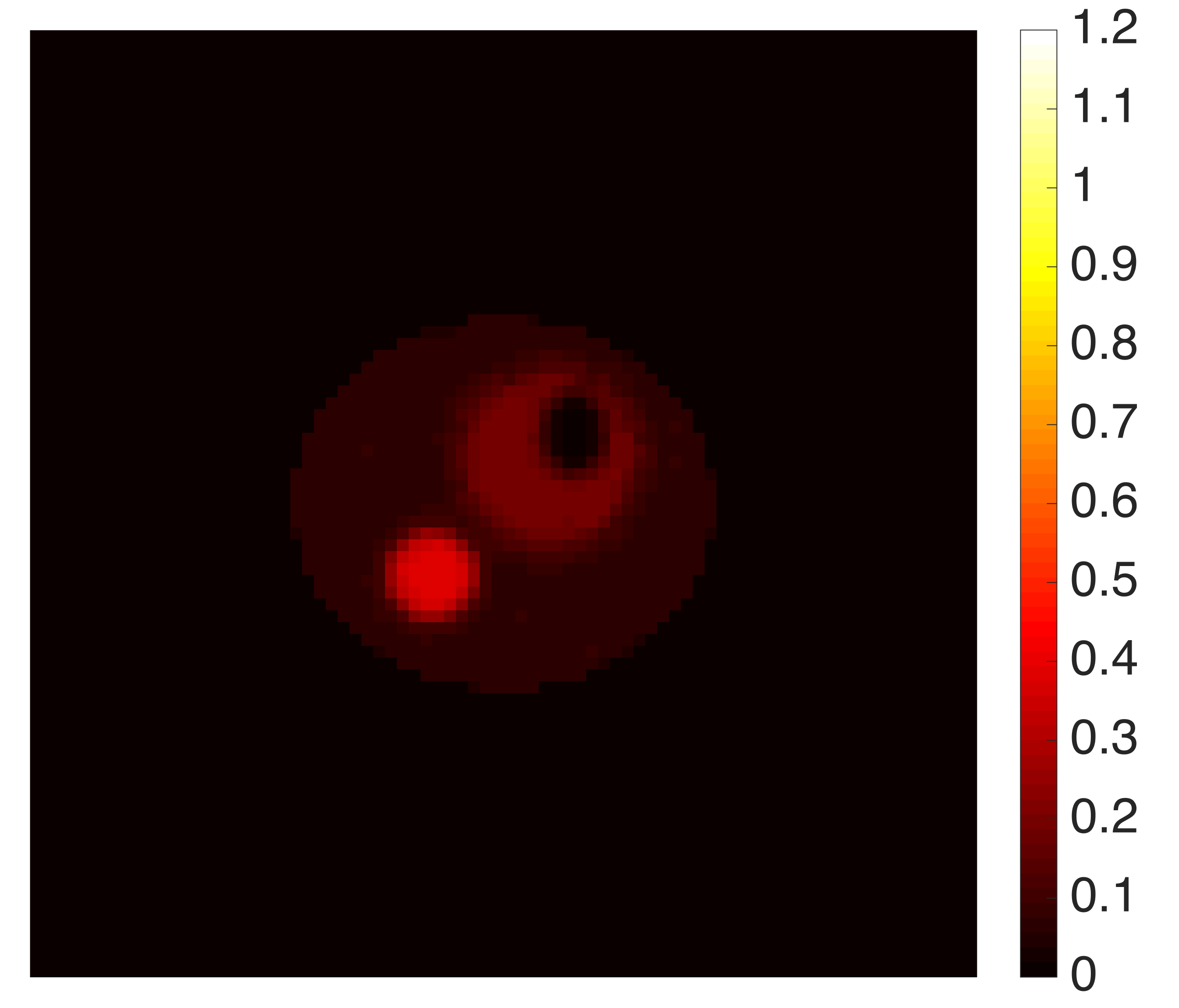}} \\
\subfigure[$\boldsymbol{K}_{fb}$ \label{fig:K1}]
{\includegraphics[width=3.8cm]{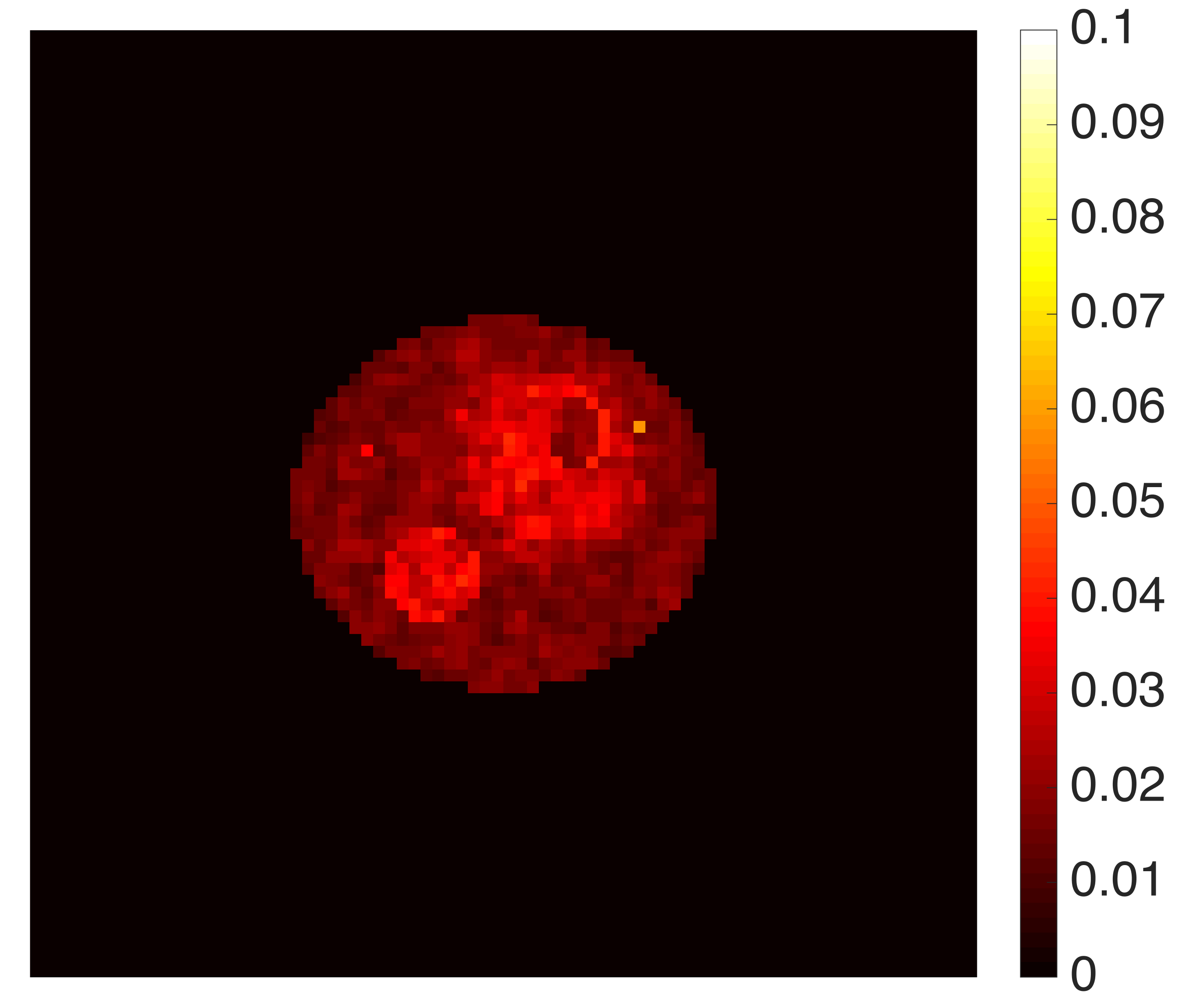}} 
\subfigure[$\boldsymbol{K}_{bf}$ \label{fig:K2}]
{\includegraphics[width=3.8cm]{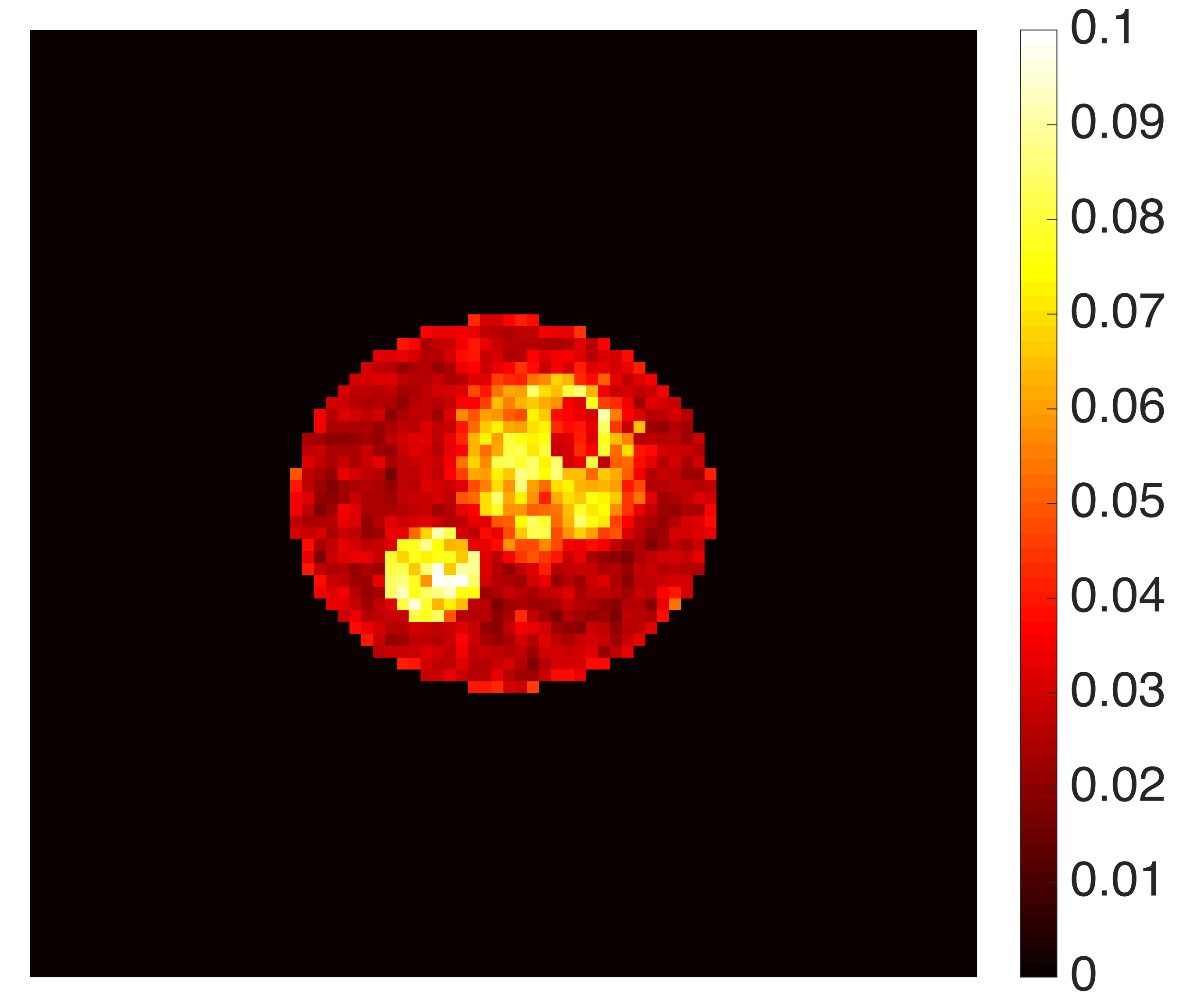}} 
\subfigure[$\boldsymbol{K}_{mf}$ \label{fig:K3}]
{\includegraphics[width=3.8cm]{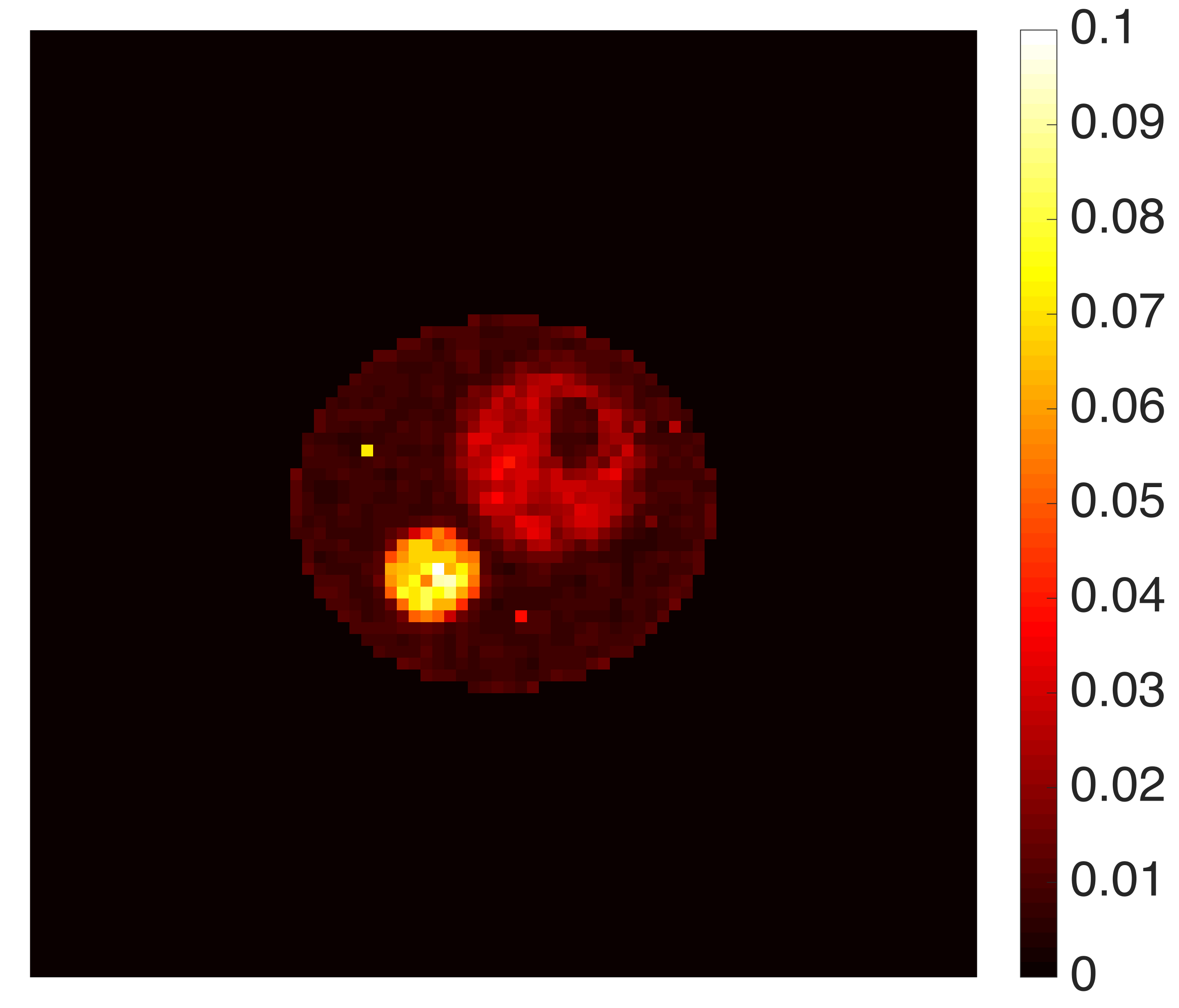}} 
\subfigure[$\boldsymbol{K}_{fm}$ \label{fig:K4}]
{\includegraphics[width=3.8cm]{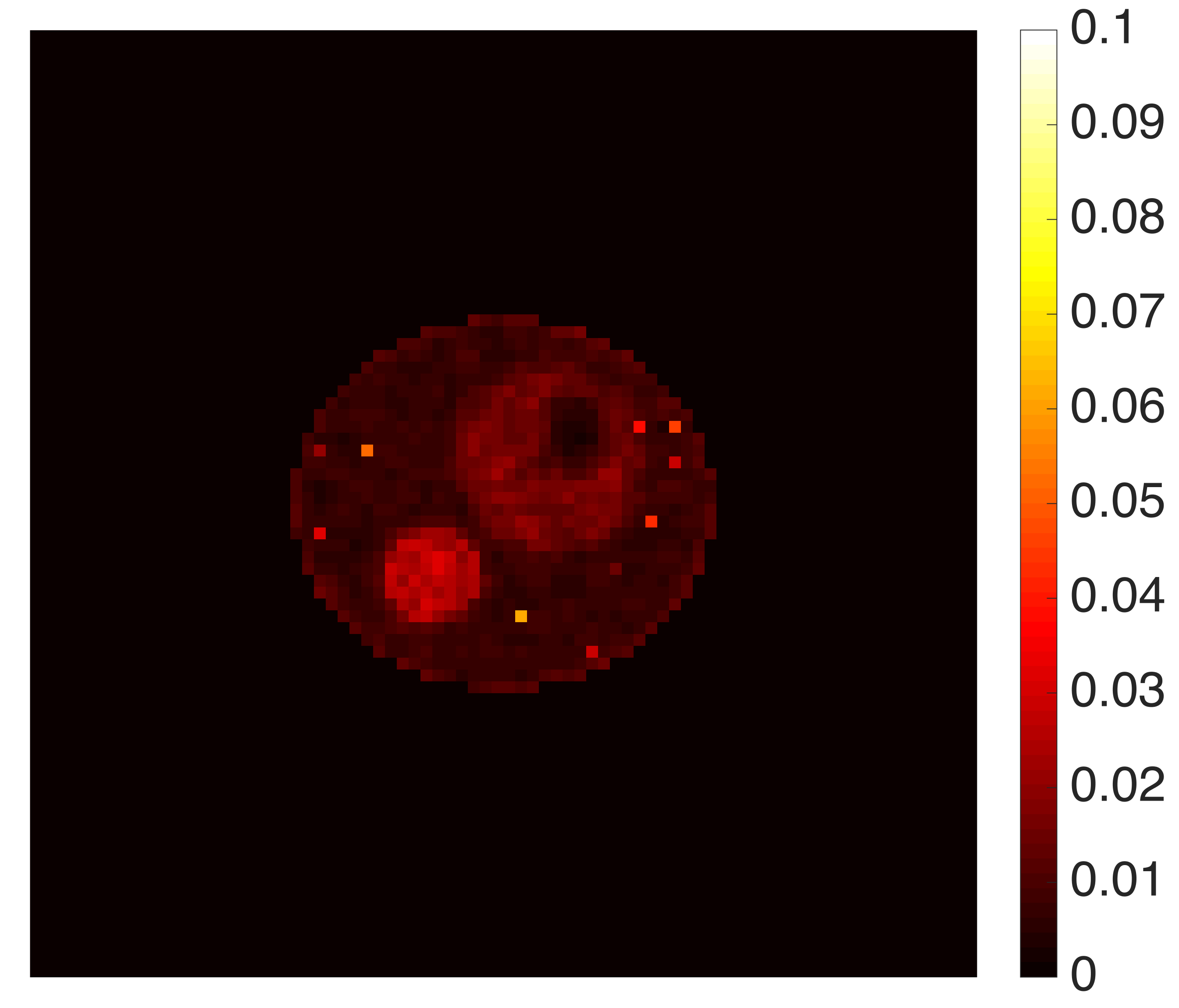}}
\caption{Mean images (first row) and standard deviation images (second row) of $\boldsymbol{K}_{fb},\boldsymbol{K}_{bf},\boldsymbol{K}_{mf},\boldsymbol{K}_{fm}$, computed over the fifty reconstructions.}
\label{fig:parametric_images_mean_std}
\end{figure}

Once the entire set of kinetic parameters $k_{fb},k_{bf},k_{mf},k_{fm}$ for each pixel were retrieved, we built up the parametric images $\boldsymbol{K}_{fb},\boldsymbol{K}_{bf},\boldsymbol{K}_{mf},\boldsymbol{K}_{fm}$.
\figurename~\ref{fig:parametric_images_mean_std} shows the mean images (first row) and the standard deviation images (second row), computed over the fifty reconstructions. The mean images provide a reliable approximation of the ground truth parametric images, demonstrating the consistency of the parametric inversion procedure. Notice that the artifacts occurring at the edges of the homogeneous regions, observable especially around the first and second regions, are consequences of the application of the Gaussian filter. The standard deviation images keep systematically small values, proving that the iterative reconstruction scheme is numerically stable with respect to noise.
\tablename~\ref{tab:kinetics_rec} reports the mean and the standard deviation of the kinetic parameters over the four homogeneous regions. Comparison between the ground truth values of \tablename~\ref{tab:kinetics} and the reconstructed values of \tablename~\ref{tab:kinetics_rec} clearly shows the reliability of our approach.

From the computational viewpoint, the parametric reconstruction takes almost 45 minutes. Please note that the algorithm was implemented in the Matlab programming environment and the algorithm was executed on a computer with a processor Intel core i5. Despite that, for a single pixel, the Gauss--Newton iterative scheme requires about 5-10 iterations before it converges and the operations carried out in a single iteration for computing the Newton step--size are not computationally demanding (the matrices in the game have small size). Therefore, the high computational cost of the method depends only from the application of the reduction scheme on a dense set of pixels. Nevertheless, the computational complexity of our parametric imaging method is consistent with standard parametric methods.

\begin{table}[htb] 
\caption{\label{tab:kinetics_rec} Mean and standard deviation of the kinetic parameters $k_{fb},k_{bf},k_{mf},k_{fm}$ (min$^{-1}$),
for each one of the four homogeneous region, computed over the fifty reconstructions.}
\begin{indented}
\item[] \begin{tabular}{@{}ccccc}
\br
& $k_{fb}$ & $k_{bf}$ & $k_{mf}$ & $k_{fm}$ \\
\mr
region 1 & $0.75 \pm 0.02$ & $0.59 \pm 0.03$ & $0.07 \pm 0.01$ & $0.07 \pm 0.01$ \\
region 2 & $0.93 \pm 0.03$ & $0.87 \pm 0.06$ & $0.16 \pm 0.02$ & $0.16 \pm 0.01$ \\                             
region 3 & $1.04 \pm 0.03$ & $0.79 \pm 0.08$ & $0.33 \pm 0.06$ & $0.31 \pm 0.02$ \\                             
region 4 & $0.58 \pm 0.02$ & $0.51 \pm 0.03$ & $0.08 \pm 0.01$ & $0.02 \pm 0.01$ \\                             
\br
\end{tabular}
\end{indented}
\end{table} 

\subsection{Real data}
We now test our parametric imaging method on real FDG--PET experiments in the case of the three--compartment non--catenary model representing the renal physiology, described in Subsection~\ref{subsec:3C}.   

We analyzed FDG--PET real data of murine models obtained by means of the dedicated microPET system (Albira, Carestream Health, Genova) \cite{Albira} currently operational at our lab. 
Following the experimental protocol for FDG--PET experiments, utilized during a study on the metabolic effects of metformin \cite{protocol}, the animals were studied after a fasting period of six hours to ensure a steady state of substrate and hormones governing glucose metabolism. Then, the animals were properly anesthetized and positioned on the bed of the microPET system whose two--ring configuration covers the whole animal body in a single bed position. A dose of 3 to 4 MBq of FDG was injected through the tail vein, soon after the start of a dynamic list mode acquisition lasting 40 min. The acquisition was reconstructed using the following framing rate: $10 \times 15$s, $1 \times 22$s, $4 \times 30$s, $5 \times 60$s, $2 \times 150$s and $5 \times 300$s. The dynamic PET images of tracer concentration (kBq/ml) were reconstructed using a Maximum Likelihood Expectation Maximization (MLEM) method \cite{Shepp}. The complete dataset is composed by 100 images of $80 \times 80$ pixels, each one reproducing an axial section, by the total number of time points of the experiment. 
For this test, we considered a mouse in a control (CTR) condition and a mouse in a starved (STS) condition (food deprivation, with free access to water, for 48 hours). We focused on the analysis of the renal physiology and selected a single PET slice containing an axial section of the right kidney, the same slice for both animals. 
The entire FDG kinetic process was initialized by the arterial IF. We are aware that the determination of IF is a challenging task in the case of mice. To accomplish it, for each animal model we have first viewed the tracer pass in cine mode. Then, in a frame where the left ventricle was particularly visible, we have drawn a ROI in the aortic arc and maintained it for all time points. For both analysis, the blood volume fraction was assumed to be equal to $V_b = 0.2$, a typical value for the kidney of the mouse \cite{Garbarino_kidney}.

\begin{figure}[htb]
\centering
\subfigure[$\boldsymbol{K}_{fa}$ \label{fig:K1_ctr}]
{\includegraphics[width=3.8cm]{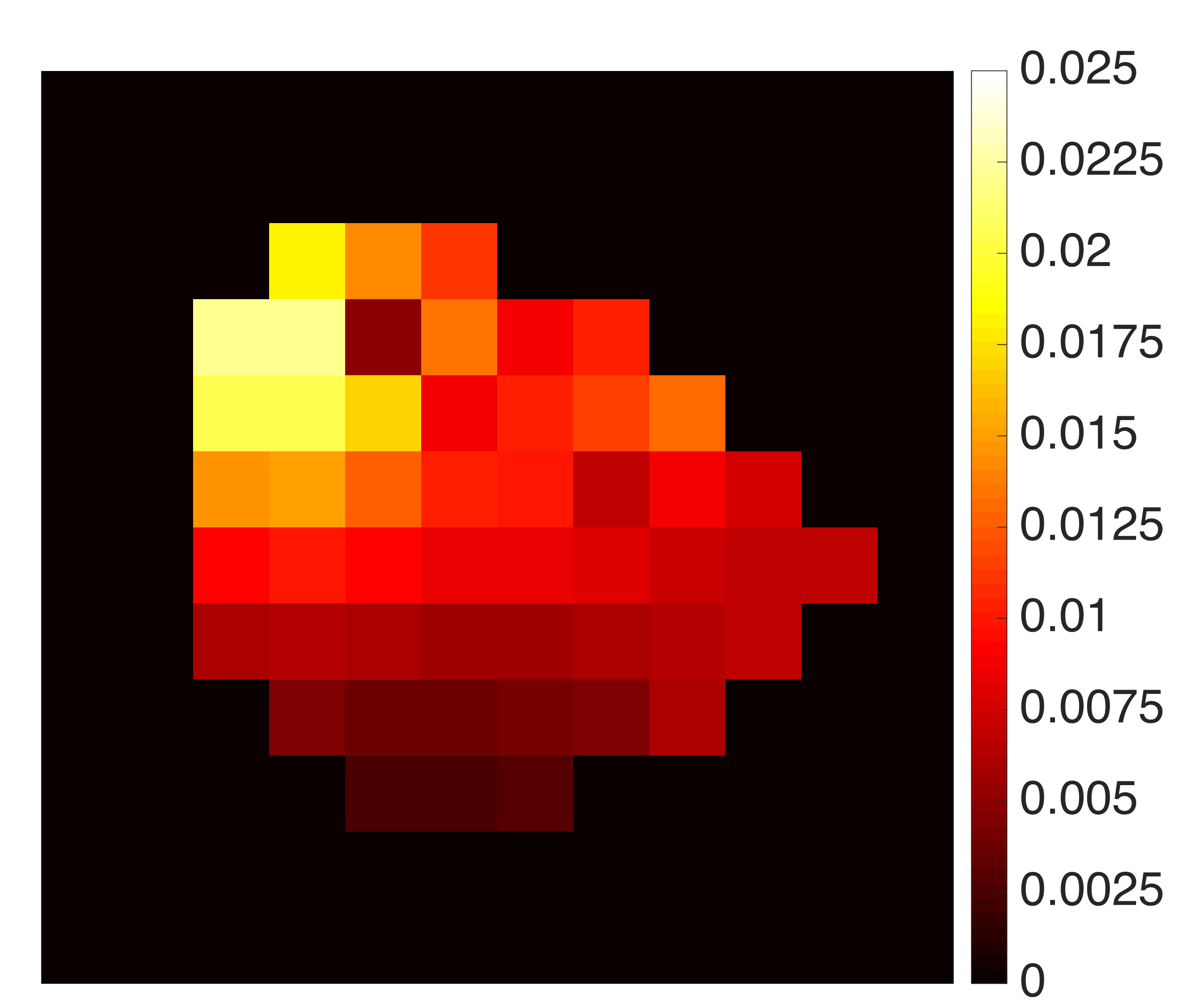}} 
\subfigure[$\boldsymbol{K}_{ma}$ \label{fig:K2_ctr}]
{\includegraphics[width=3.8cm]{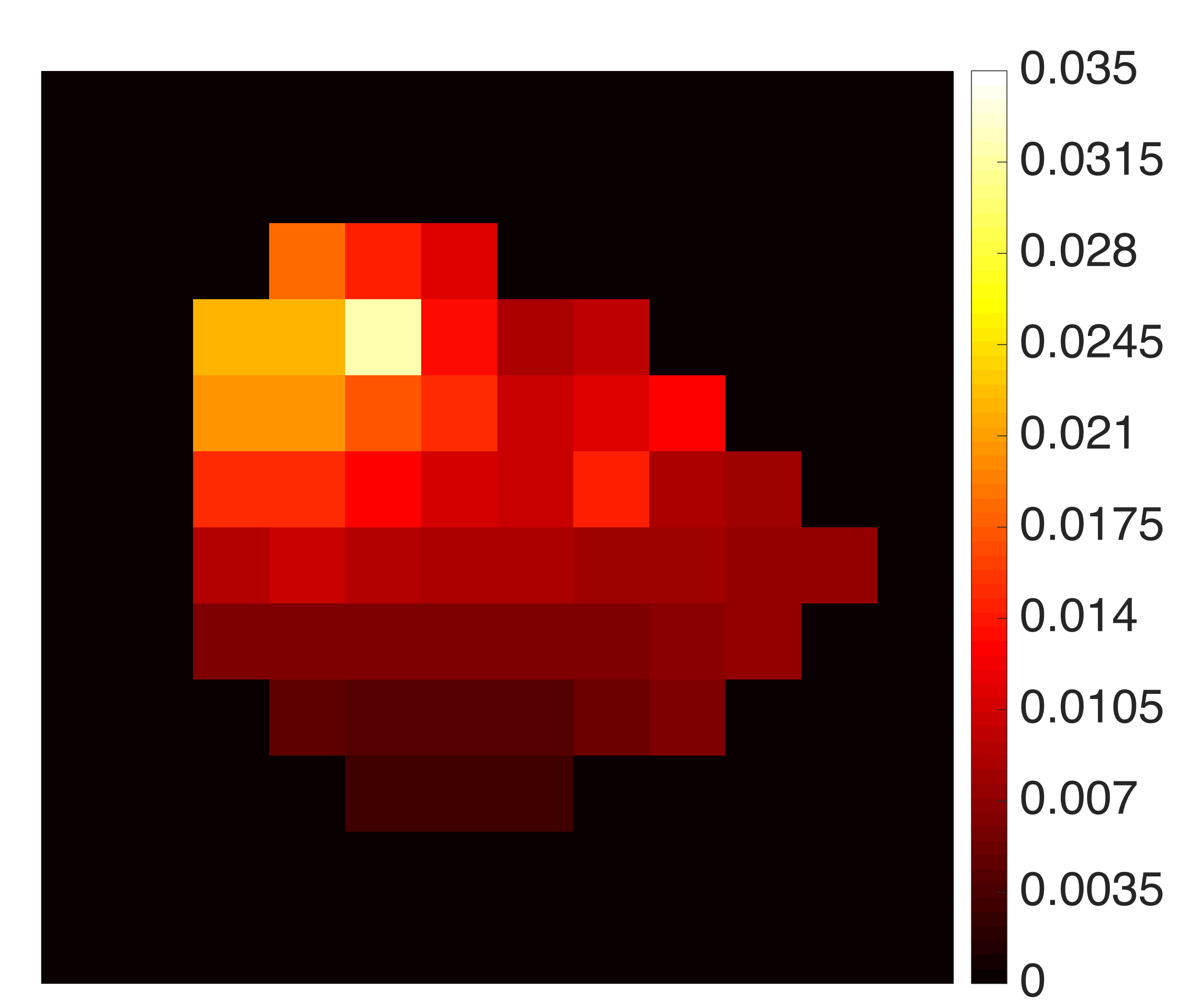}} 
\subfigure[$\boldsymbol{K}_{af}$ \label{fig:K3_ctr}]
{\includegraphics[width=3.8cm]{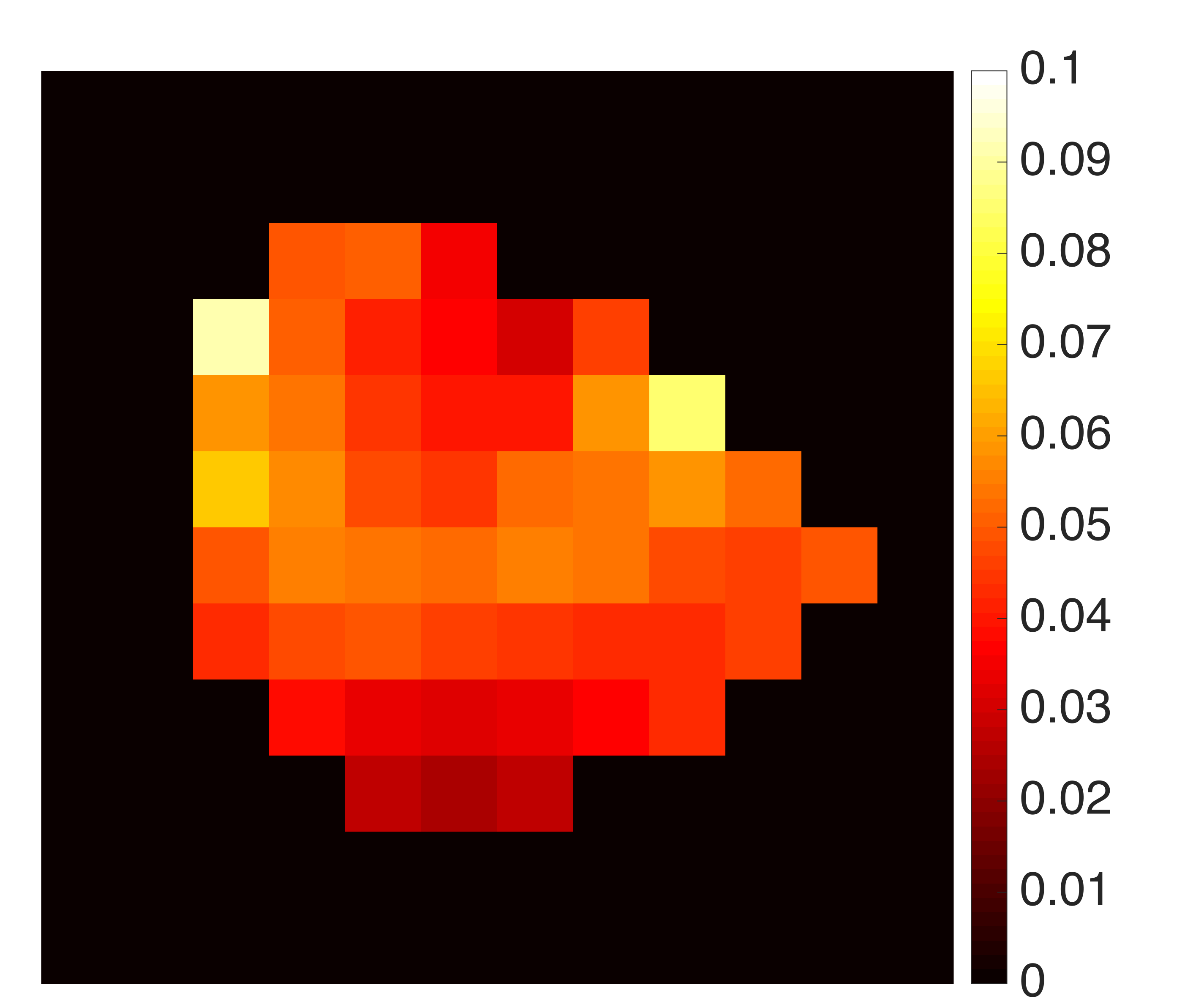}} \\
\subfigure[$\boldsymbol{K}_{fa}$ \label{fig:K1_sts}]
{\includegraphics[width=3.8cm]{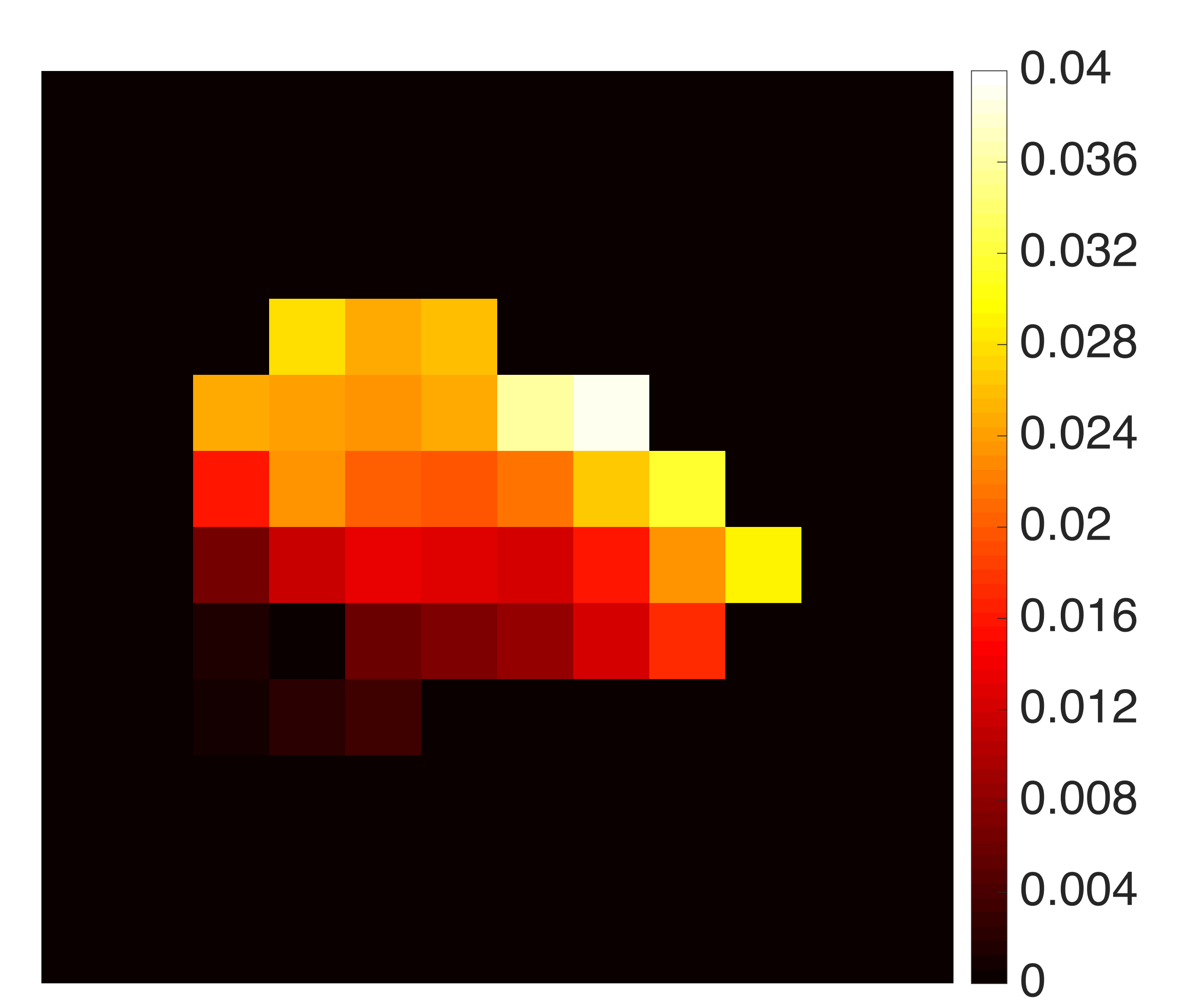}} 
\subfigure[$\boldsymbol{K}_{ma}$ \label{fig:K2_sts}]
{\includegraphics[width=3.8cm]{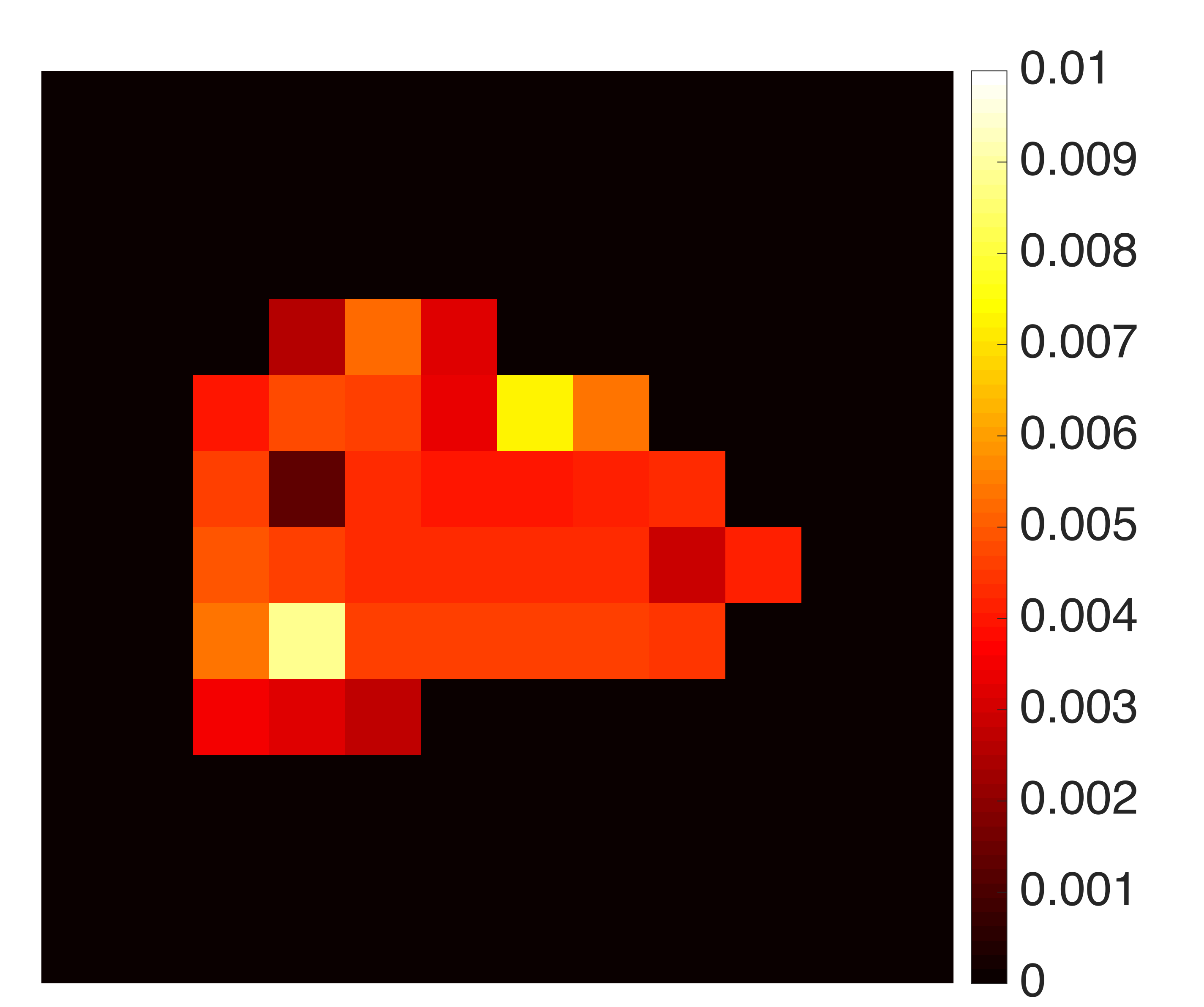}} 
\subfigure[$\boldsymbol{K}_{af}$ \label{fig:K3_sts}]
{\includegraphics[width=3.8cm]{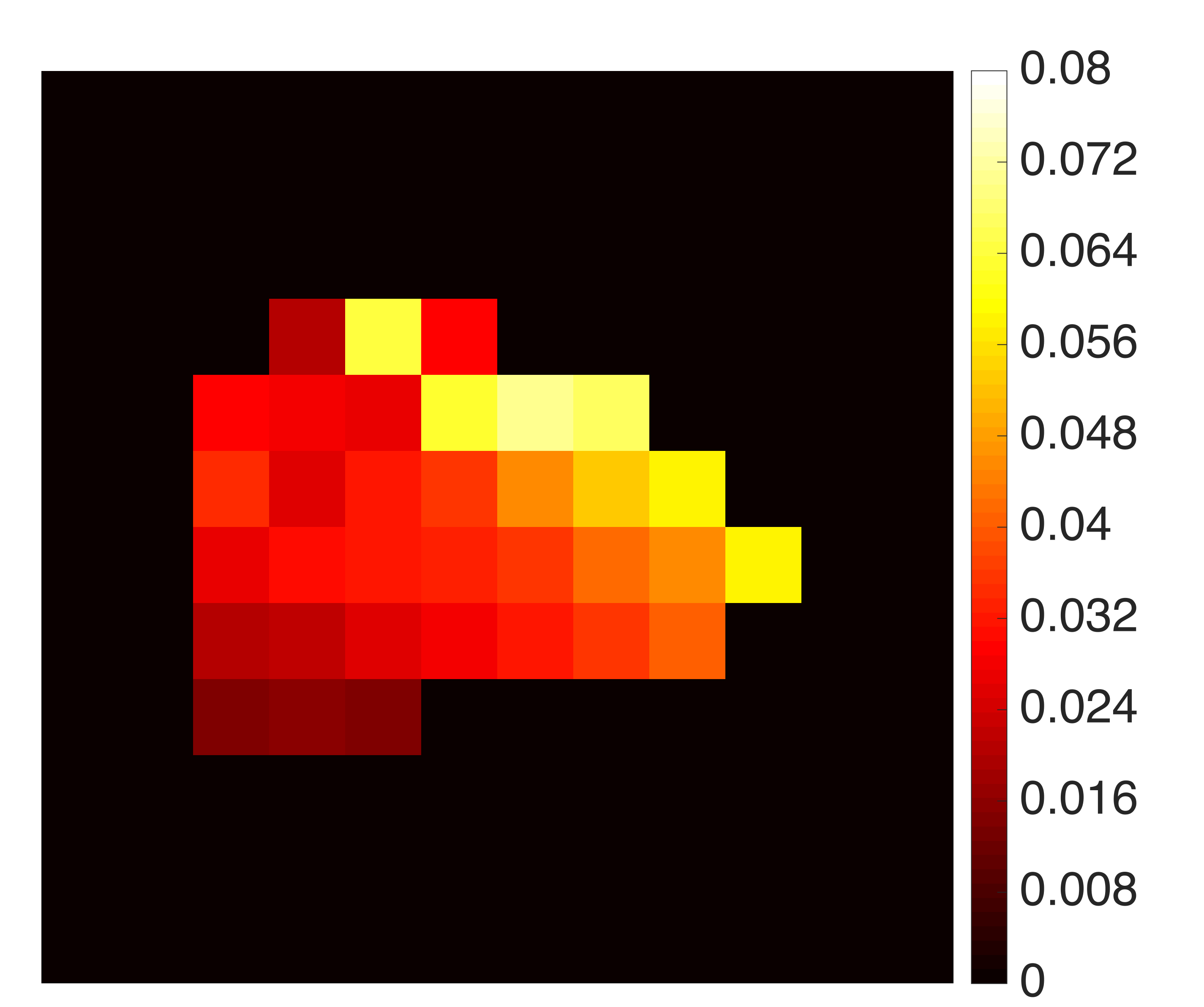}} 
\caption{Parametric images $\boldsymbol{K}_{fa},\boldsymbol{K}_{ma},\boldsymbol{K}_{af}$: first row for the CTR mouse, second row for the STS mouse.}
\label{fig:parametric_images_kidney1}
\end{figure}

\begin{figure}[htb]
\centering
\subfigure[$\boldsymbol{K}_{mf}$ \label{fig:K4_ctr}]
{\includegraphics[width=3.8cm]{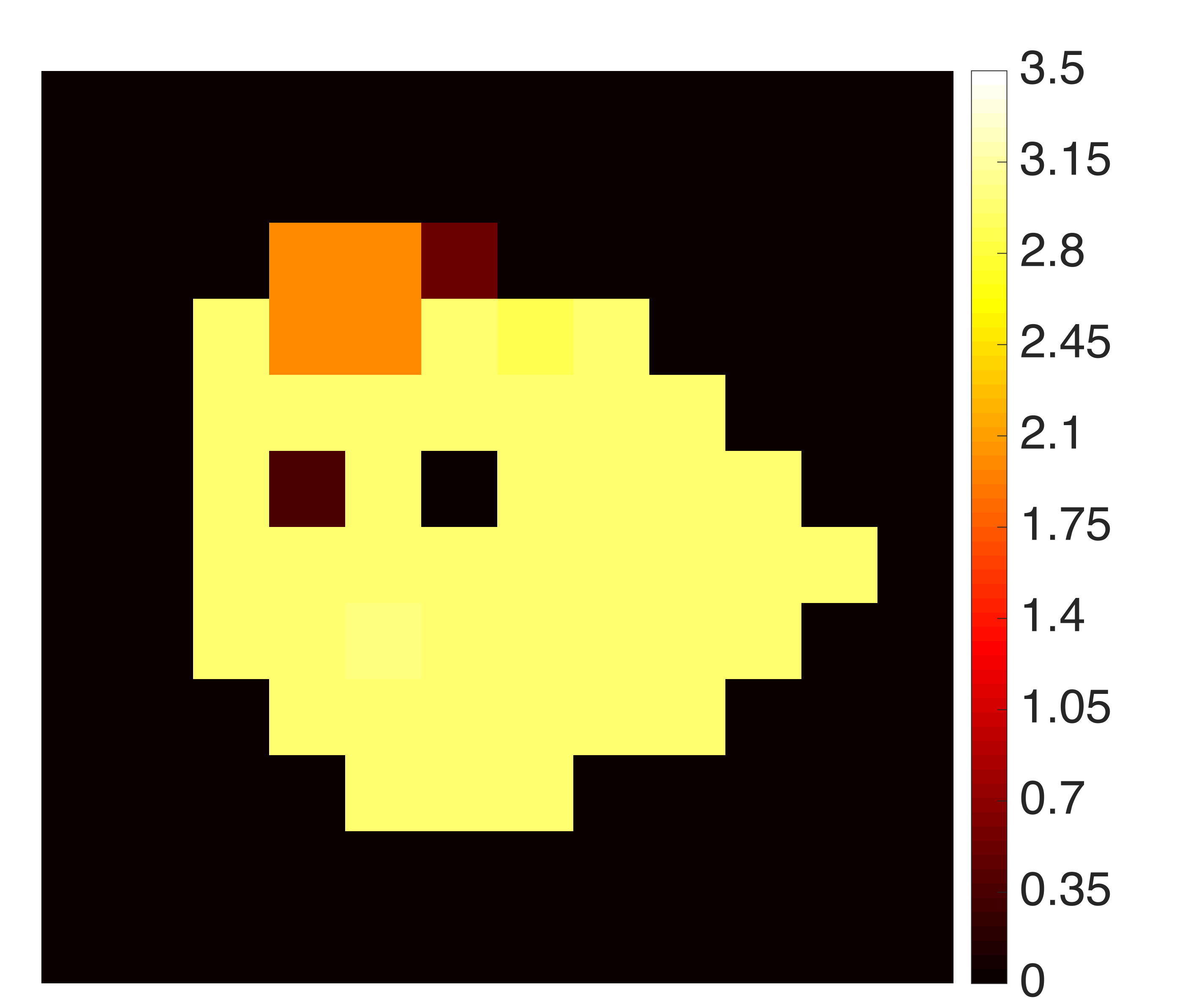}} 
\subfigure[$\boldsymbol{K}_{fm}$ \label{fig:K5_ctr}]
{\includegraphics[width=3.8cm]{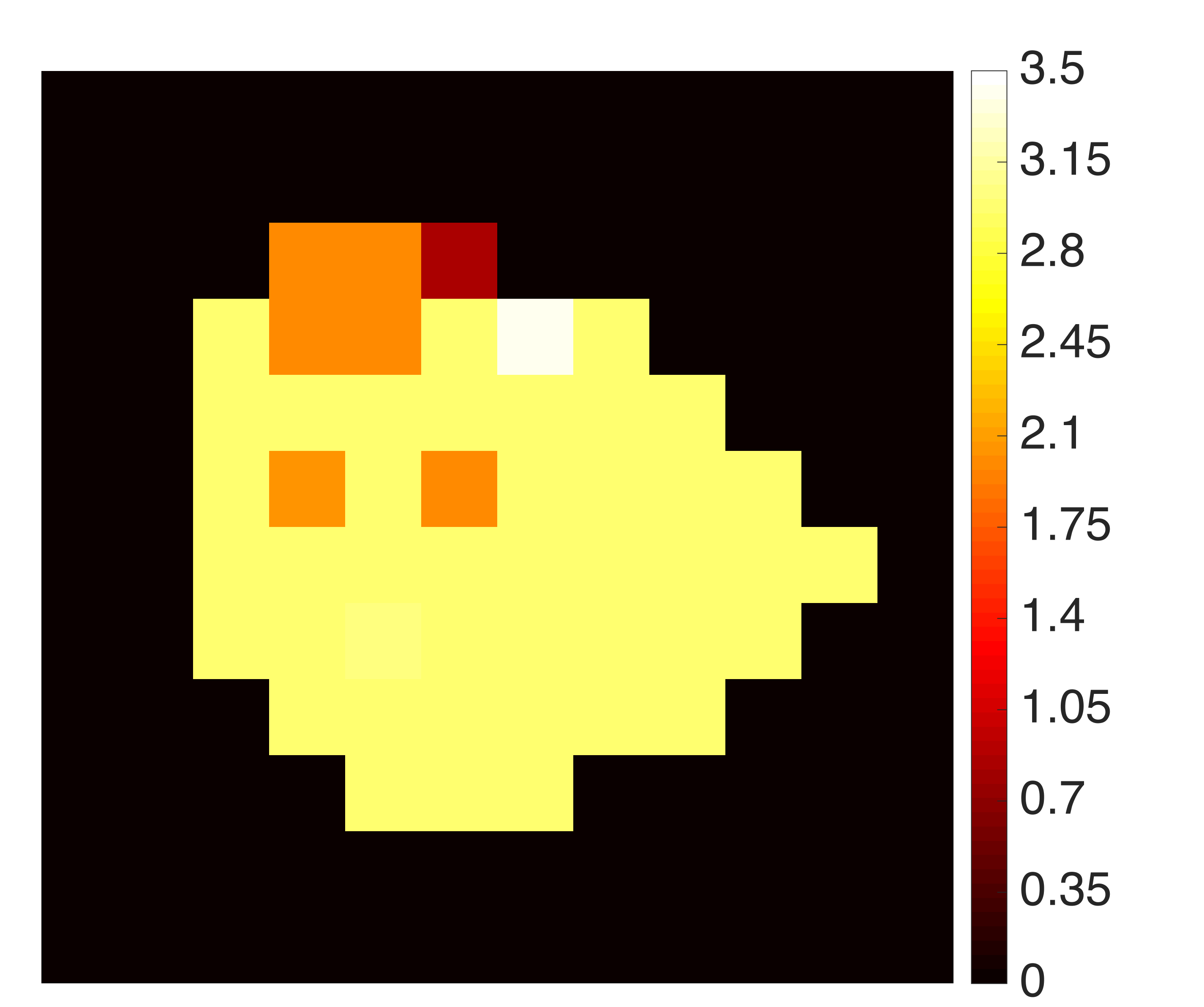}} 
\subfigure[$\boldsymbol{K}_{tm}$ \label{fig:K6_ctr}]
{\includegraphics[width=3.8cm]{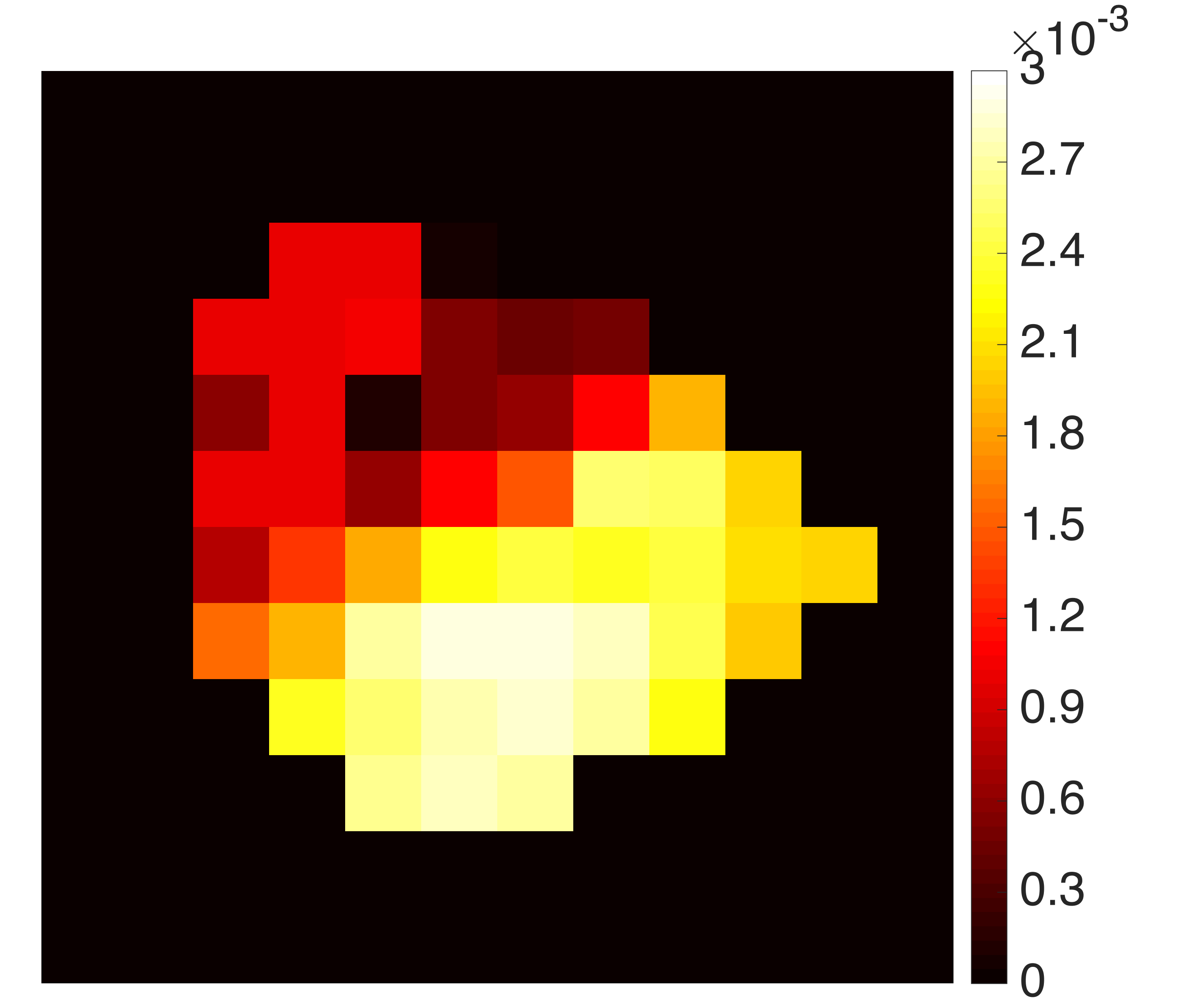}} 
\subfigure[$\boldsymbol{K}_{ut}$ \label{fig:K7_ctr}]
{\includegraphics[width=3.8cm]{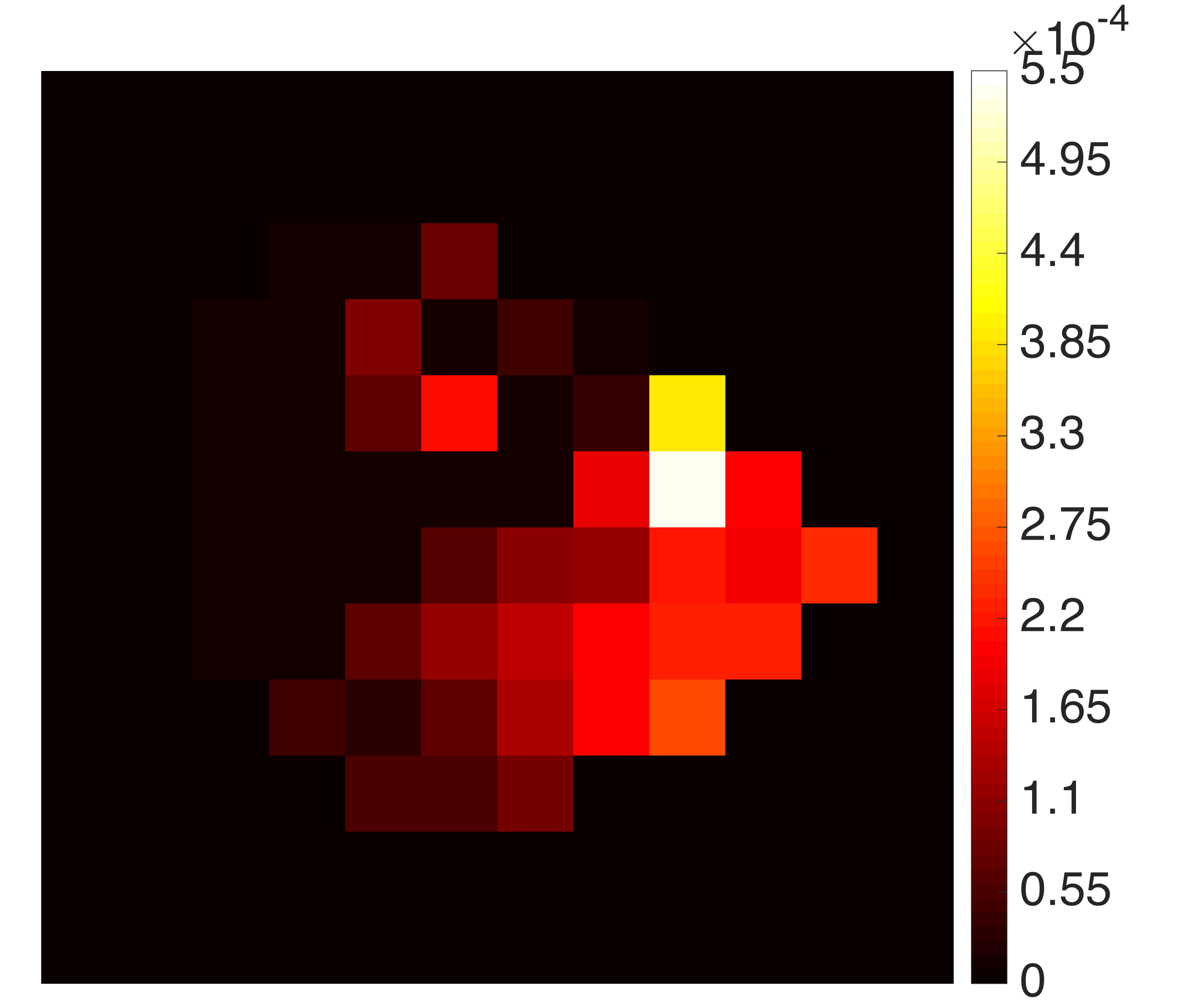}} \\
\subfigure[$\boldsymbol{K}_{mf}$ \label{fig:K4_sts}]
{\includegraphics[width=3.8cm]{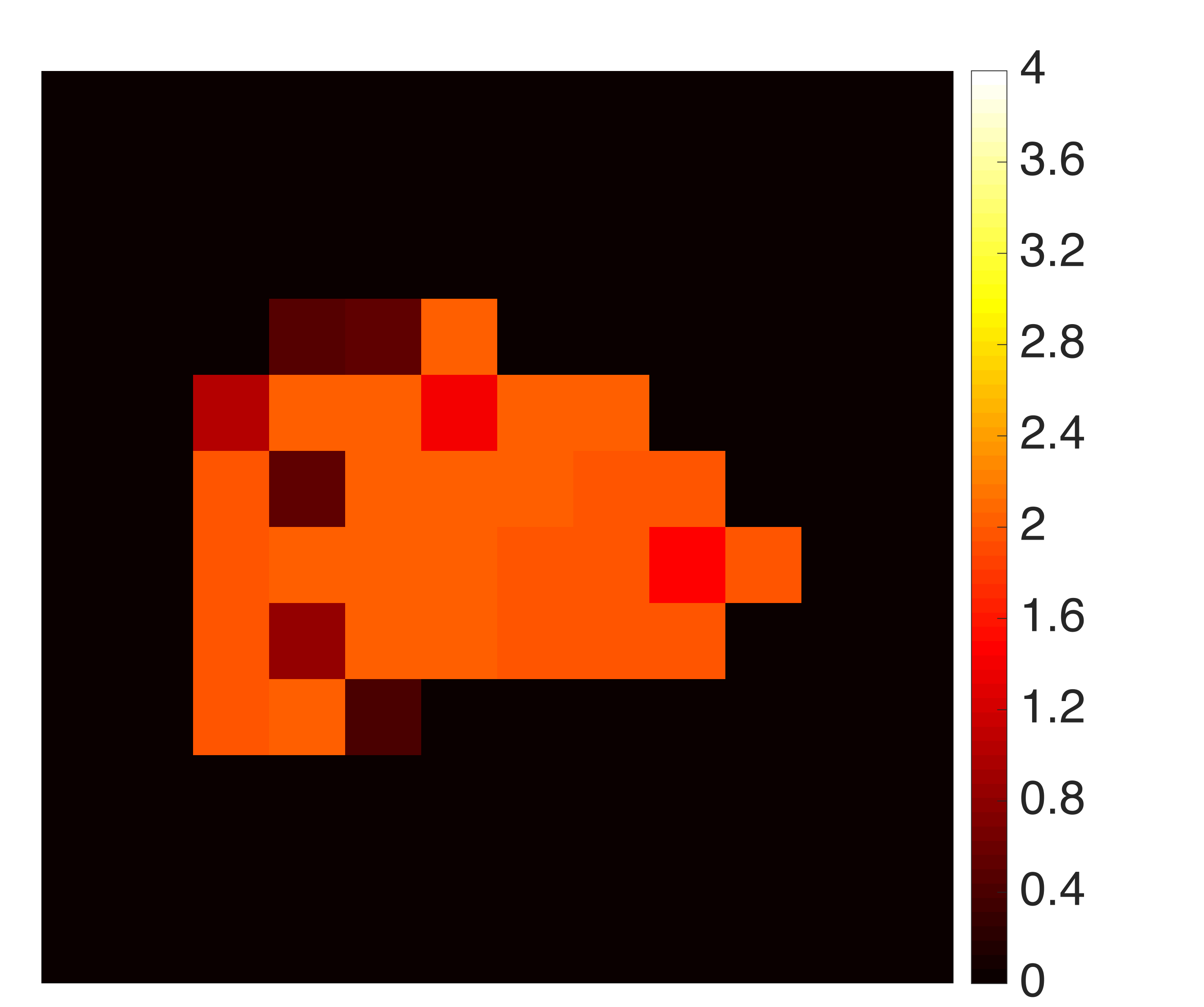}} 
\subfigure[$\boldsymbol{K}_{fm}$ \label{fig:K5_sts}]
{\includegraphics[width=3.8cm]{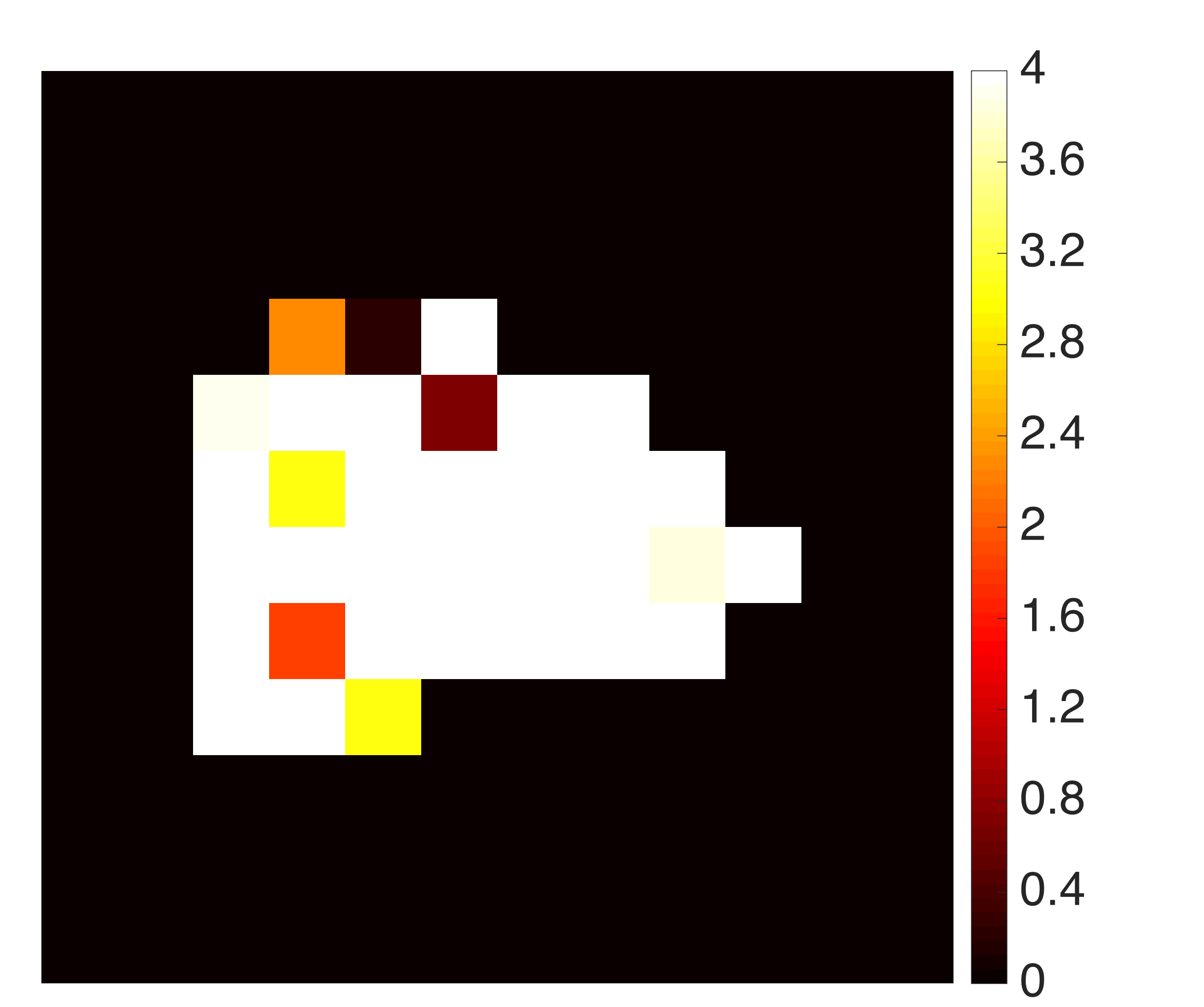}} 
\subfigure[$\boldsymbol{K}_{tm}$ \label{fig:K6_sts}]
{\includegraphics[width=3.8cm]{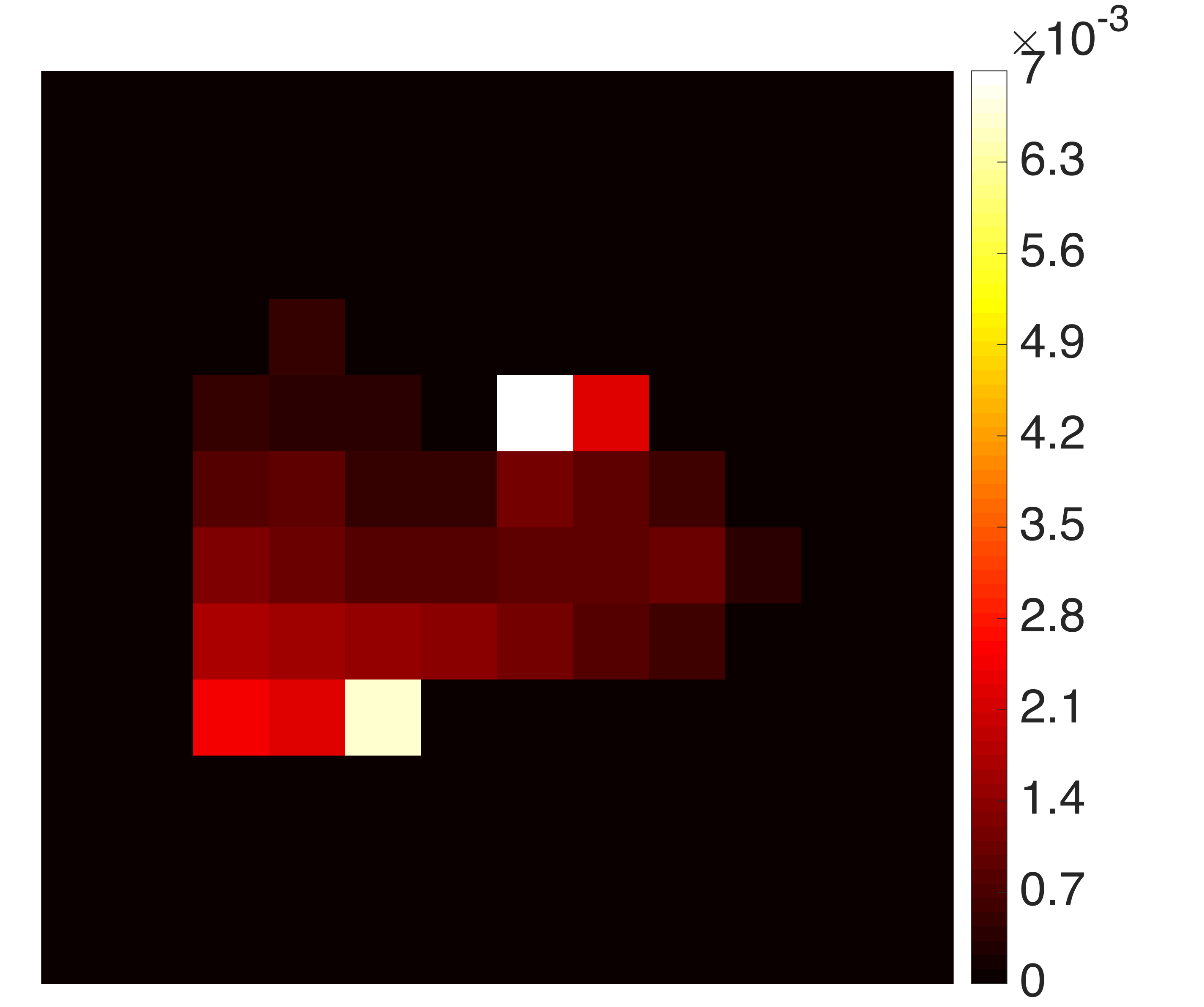}} 
\subfigure[$\boldsymbol{K}_{ut}$ \label{fig:K7_sts}]
{\includegraphics[width=3.8cm]{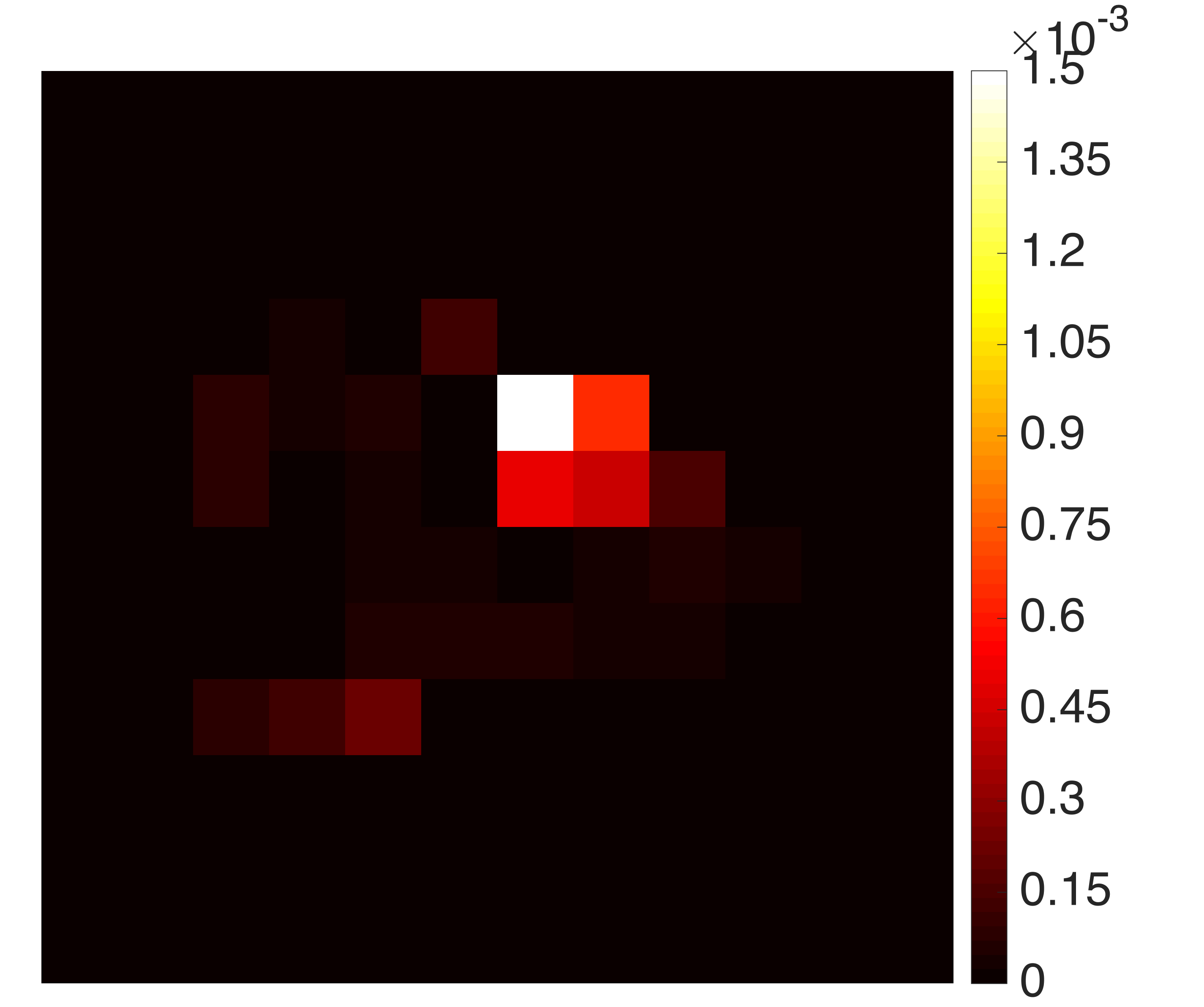}} 
\caption{Parametric images $\boldsymbol{K}_{mf},\boldsymbol{K}_{fm},\boldsymbol{K}_{tm},\boldsymbol{K}_{ut}$: first row for the CTR mouse, second row for the STS mouse.}
\label{fig:parametric_images_kidney2}
\end{figure}

We applied our imaging method on the selected dynamic PET slice of the CTR mouse and of the STS mouse. More specifically, we smoothed the data by means of a Gaussian filter of standard deviation $\sigma = 1$ and size $3 \times 3$; we selected the ROI within the axial section of the kidney through the image segmentation process determining the minimum value of activity recorded inside the organ, and reconstructed the kinetic parameters $k_{fa}, k_{ma}, k_{af}, k_{mf}, k_{fm}, k_{tm}, k_{ut}$ of the three--compartment system for each pixel by means of the regularized Gauss--Newton iterative algorithm. The initial guesses were randomly selected in the interval $(0,1)$ and the regularization parameter was chosen by the GCV method.
The reconstructed parametric images $\boldsymbol{K}_{fa},\boldsymbol{K}_{ma},\boldsymbol{K}_{af},\boldsymbol{K}_{mf},\boldsymbol{K}_{fm},\boldsymbol{K}_{tm},\boldsymbol{K}_{ut}$ for the renal compartmental model are presented in \figurename~\ref{fig:parametric_images_kidney1} and \figurename~\ref{fig:parametric_images_kidney2}: in each figure, the first row shows the parametric images for the CTR mouse, the second row for the STS mouse.
     
All parametric images, of both the CTR mouse and the STS mouse, show parameters' values that vary quite largely from pixel to pixel, bringing out the lack of homogeneity of the renal tissue. Indeed, the parametric images point out the different structures composing the kidney and characterizing the distinct functions of the organ \cite{kid_phys}. 
This is consistent with the architecture of the renal compartmental model we have designed (Subsection~\ref{subsec:3C}).
The higher activity of the parameters is located in a specific part of the outer portion of the axial section of the kidney, which is attributable to the renal cortex, in which most of the renal processes are carried out. Moreover, for both the CTR and STS conditions, we can observe that the parametric images $\boldsymbol{K}_{tm}$ and $\boldsymbol{K}_{ut}$ linked to the tubule compartment have a very similar distribution while the physiologically sound condition $k_{tm} \simeq 10^{2} k_{ut}$ is maintained pixel--wise (without any constraint in the inversion procedure).
Instead, the fundamental difference between the CTR and the STS parametric images relies on the numerical values of the parameters. 
In particular, the values of the exchange coefficients associated with the blood input, that are $k_{fa}, k_{ma}, k_{af}$, exhibit a different behavior with respect to the two conditions. From the CTR mouse to the STS mouse: the input parameter from the blood $k_{fa}$ increases (almost duplicates), the filtration process described by $k_{ma}$ decreases (is almost reduced by a factor of three), and the output parameter to the blood $k_{af}$ remains approximately equal.
This trend reflects the response of the kidney to the different physiological conditions of the two mice analyzed, coherently with what already observed \cite{Garbarino_kidney}.
Finally, we notice that the parameters linked to the FDG metabolization process, $k_{mf}$ and $k_{fm}$, and the parameters representing the reabsorption and excretion processes, $k_{tm}$ and $k_{ut}$ respectively, remain basically unchanged in the two conditions.

\section{Conclusions}
FDG--PET imaging allows the observation of metabolic processes related to glucose consumption inside a living organism. In order to improve the quality of information achievable from PET images, one solution is to develop parametric imaging methods capable of showing the tracer metabolism pixel--wise. Starting from the design of compartmental models suitable to describe the tracer kinetics in a predefined physiological system, parametric imaging techniques process dynamic PET images and estimate the spatial distribution of the exchange coefficients identified by the model.

In this paper we have shown a novel parametric imaging tool, which integrates pre--processing methods and optimization of non--linear inverse problems to realize an automatic pipeline for the reduction of multi--compartment models. The novelties of this approach are at a mathematical, computational, and technological level. In fact, from a mathematical viewpoint, we showed that model identifiability occurs also in the case of more complex compartmental systems like the non--catenary one involving three compartments and nicely mimicking the renal physiology. From a computational viewpoint, we showed that the use of a regularized Gauss--Newton method for the reduction of the compartmental models is effective also in the case of a pixel--wise analysis. Finally, from a technological viewpoint, we provided a simple pipeline, which is able to realize the complete processing workflow, starting from dynamical nuclear medicine data up to images of all kinetic parameters. The main advantage of this approach is in its notable degree of generality, since in principle it may be applied to models made of several compartments. Further, differently to typical linear parametric imaging methods, this algorithm provides maps of all model parameters, and differently to direct imaging methods, it may reconstruct a large set of kinetic parameters, and account for different models in the same image. We focused on the standard two--compartment catenary model and a three--compartment non--catenary model representing the renal system. Interestingly, in this latest application, the results obtained from experimental observations showed that the approach works properly without the need of using data from the bladder, which implies that we were able to obtain reliable parametric reconstructions using just measurements from pixels in the kidneys. We also remark that the procedure can be applied also for a large variety of compartmental models, more general and complicated, provided that the numerical complexity of the model is taken into account.
We validated our approach against synthetic data and then illustrated it on FDG--PET data of murine models. We analyzed dynamic FDG--PET images of a selected axial section of the right kidney for a control mouse and a starved mouse, proving that the reconstructed renal parametric images are qualitatively effective in the description of the local FDG metabolism.

Further developments of this approach to parametric images of nuclear medicine data are concerned with several issues. From the computational viewpoint, we aim at reducing the computational burden of this approach by means of \emph{ad hoc} designed implementations. This will allow the implementation of a technological pipeline able to deal with four-dimensional data. Moreover, we are working at the extension of this approach to the direct imaging problem, whereby the input data are the PET count sinograms and not the reconstructed PET images. Finally, from the clinal viewpoint, we are going to apply this approach against a notable quantity of data acquired according to paradigms designed to investigate the role of metformin in glucose metabolism \cite{protocol}.  


\end{document}